\input amstex
\documentstyle{amsppt}
\input boxedeps.tex
\SetepsfEPSFSpecial

\HideDisplacementBoxes
\hsize=5.5in
\vsize=7.4in
\magnification=1200
\def\e{\epsilon}
\def\p{^\prime}
\def\pp{^{\prime \prime}}
\def\i{^{-1}}

\def\eg{\eurm g}
\def\ep{\eurm p}
\def\eS{\eurm S}
\def\eM{\eurm M}
\def\en{\eurm n}
\def\eB{\eurm B}

\def\eb{\eurm b}
\def\a{\alpha}
\def\g{\gamma}
\def\s{\sigma}
\def\b{\beta}
\def\d{\delta}
\def\out{{\text{out}}}
\def\id{{\text{id}}}
\parindent=18pt
\loadeurm

\document
\baselineskip=18pt
\NoBlackBoxes
\def\leaderfill{\leaders\hbox to 1em{\hss.\hss}\hfill}

{\centerline{\bf{ON COMPLEXITY OF THE WORD PROBLEM}}}
{\centerline{\bf{IN BRAID GROUPS AND MAPPING CLASS GROUPS}}}

\

{\centerline{{Hessam Hamidi-Tehrani}
\footnote{Department of Mathematics, University of California at
Santa Barbara, CA 93106
E-mail: hessam\@math.ucsb.edu.
http://www.math.ucsb.edu/$\sim$hessam.}}}

\

{\centerline{{September 25, 1998}}}

\

\topmatter
\abstract
We prove that the word problem in the mapping class group of the once-punctured surface 
of genus $\eurm g$ has complexity 
$O(|w|^2 \eurm g)$ for $|w| \ge \log(\eg)$ where $|w|$ is the length of
the word in a (standard) set of generators. The corresponding bound 
 in the
 case of the closed surface is $O(|w|^2 \eg^2)$.  
We also carry out the same methods 
for the braid groups, and show that this gives a bound which improves the best known bound in this case; namely, the complexity of the word problem in
the $\en-$braid group is $O(|w|^2 \en)$, for $|w| \ge \log \en$.
 We state a similar result for mapping class groups of surfaces with
several punctures.

\


\


\endabstract
\endtopmatter

\head \S 0. Introduction  \endhead

A group $G$ is said to have a {\it{solvable word problem}} if
there is a finite generating set $S$ for $G$ such that there is an algorithm
to decide if a given word $w$ in $S$ represents the identity element in $G$.
The word problem is said to have {\it{complexity}} $O(f(|w|))$
if there exist such an algorithm which takes $\le k f(|w|)$ steps
on a Turing Machine (TM) to produce a ``yes'' or a ``no'', for a word $w$ of length $|w|$ 
where $k$ is a constant (see Appendix for more on complexity and Turing Machine).
The conjugacy problem is defined similarly, but the objective is to 
decide if two given words are conjugate in the group $G$.

Sometimes one has to deal with sequence of groups $G_n$ depending on an
integer parameter  $n$ (say mapping class groups of closed surfaces which
is parameterized by genus), and one can pose the question of how the
complexity of a problem grows as  $n$
 becomes larger. This is a crucial issue in implementation of a uniform algorithm, because the parameter becomes an input. In this case we say the word problem has {\it{uniform complexity}} $O(f(|w|,n))$ for the groups $G_n$ if 
there exists some finite set of generators for each $G_n$ such that for a word $w$ in generators of $G_n$ of length $|w|$, it takes a Turing Machine
$\le k f(|w|, n)$ steps to determine if $w=1$. 

The word problem and conjugacy problem in the mapping class
group have been known to be solvable for a long time (see [G],[He],[P],[Mo1]). In recent years,
with development of the theory of automatic groups, some new ideas in
this direction have been discovered. In [E], the authors discuss an
automatic structure derived from Garside's algorithm [G] for the braid
groups. This results in an algorithm which is of uniform  complexity 
$O(|w|^2 \en \log \en)$,
 where $\en$ is the number 
of strands, and $|w|$ is the length of the braid, 
which is given as a word $w$
in the standard set of Artin generators (see (3.1)).
 Mosher [Mo2] proved that mapping class groups are
automatic, giving an algorithm for the word problem 
which is  quadratic in the word length
[Mo3], with no implication on uniform complexity. As the authors of [E] mention, it is important to have a bound
on the  uniform complexity; i.e., in terms of the genus and the number of 
punctures. 
Here we prove that the word problem in the mapping
 class group of the closed surface of genus $\eg$ has complexity
$O(|w|^2\eg^2+|w|\eg^2 \log \eg).$ The corresponding bound for a once-punctured 
surface of genus $\eg$ is $O(|w|^2 \eg+ |w| \eg \log \eg)$.

In a sense we  answer the Open Question 9.3.10 in [E],
but we do not use the automatic theory. Our methods rely on the action of 
the mapping class group on the space of curves, or measured train-tracks.
This could be related to the Open Question 9.4.5 in there as well, although we do not 
speak about conjugacy problem at all. It is an interesting question to try to
use the methods here to solve and analyze the complexity of the conjugacy
problem in the mapping class groups. In this respect the work of
Kleinberg and Menasco [KM], Masur and Minsky [MM1], [MM2] is of interest.
In particular, the authors of the latter prove that if two pseudo-Anosov maps are conjugate, then there is a conjugating element whose word length is linearly bounded by the larger of the word lengths of those elements.

Our methods apply to the braid groups $\eB_\en$ and give the complexity
$O( |w|^2 \en +|w| \en \log \en)$, which is the best known bound to date. In [BKL] the authors give a fast and practical algorithm for the
word problem in  $\eB_\en$, which works well with a ``Random Access Memory''
(RAM) machine, and has ``complexity'' $O(|w|^2 \en)$. But RAM is usually
much faster than TM (In particular, they 
assume that the braid index $\en$ can be 
encoded in one unit of memory; see Appendix), and their
 algorithm gives the same complexity as in [E],
namely $O(|w|^2 \en \log \en)$ if practiced on a TM.

Here is an outline of the rest of this paper: In $\S 1$ 
we develop the necessary notation 
for measured $\pi_1$-train-tracks and the mapping class groups.
In $\S 2$ we prove the bound on the complexity of the word problem in
once-punctured surfaces. In $\S 3$ we apply our methods
 to deduce a bound on the complexity of the word problem in the braid groups.
In $\S 4$ we develop the theory for closed surfaces; we prove the analog
to Theorem 1.5 for closed surfaces. $\S 5$ is devoted to analyze the complexity
of the word problem in closed surfaces. Finally in the appendix we
briefly address some issues about our definition of complexity.

{\proclaim{Acknowledgments}} Some of this work was part of my Ph.D. thesis
in  Columbia University. I thank everybody in the Mathematics Department,
especially my advisor Joan Birman for her graceful
support and guidance. Also, I found the referee's comments most valuable.

\

\head \S 1.  Some notation and background on train-tracks \endhead

Let $\eurm S=\eurm S^\eurm p_{\eurm g}$ be an oriented surface of
genus $\eurm g$ with $\eurm p$ fixed points, called punctures.
 Let $\eurm M=\eurm M_\eurm S=\eurm M^\eurm p_{\eurm g}$ the mapping
class group of $\eurm S$, i.e., the group $\Cal H(\eurm S)/\Cal
H_0(\eurm S)$, where $\Cal H(\eurm S)$ is the group of homeomorphisms
of $\eurm S$
 fixing the punctures pointwise, and $\Cal H_0(\eurm S) \subseteq \Cal
H(\eurm S)$ is the (normal) subgroup of the ones homotopic to identity
within $\Cal H(\eurm S)$. We denote the elements of $\eurm M$ by
$f,g,$ etc. An element of $\eurm{M}$ can be thought of as an isotopy
class of a homeomorphism (or diffeomorphism) of $\eurm S$. Sometimes
we pick a representative of the class $f$ and call it $f$ too. We
assume $\eurm S$ has a given smooth or piecewise linear structure,
depending on what suits the situation the best.

Notice that if $\eurm S\p$ is a surface with $b$ boundary components,
one can define the mapping class group $\eurm {M_{S\p}}$ of $\eurm
S\p$ by the group of isotopy classes of diffeomorphisms which fix the
boundary components pointwise. Let $\eurm S$ be obtained by shrinking
the boundary components of $\eurm S\p$ to punctures. Then we have the
short exact sequence $$1 \to \Bbb Z^b \to \eurm{M_{S\p}}
 \to \eurm{M_S} \to 1.\tag{1.1}$$

In the following we only study the surfaces $\eurm{S^p_g}$. The
corresponding information about surfaces with boundary can be obtained
 using (1.1).

{\proclaim{ Definition 1.1 \rm{(Train-track)}}} (See [PH].)  A
compact, connected subset $\tau$ of $\eurm S$ is called a {\it {train-track}}
if $\tau$ is a smooth branched 1-manifold embedded smoothly in $\eurm
S$. At each branch point $v$ (also called a switch point) there is a
well-defined tangent space.  Every connected component of $\tau -
\{$branch points$\}$ is called a branch. There is a natural partition
into two subsets for the set of
 branches $b$ coming to a switch $v$ (i.e., $v \in \bar {b} $)
depending on which direction they become tangent at the switch
point. We call these two sets $incoming$ and $outgoing$. The
particular choice does not matter. Also, there is a ``hyperbolicity
condition'' on the complement $\eurm S-\tau$: The doubles of
components of $\eurm S-\tau$ must have negative Euler
characteristic. Notice that the double of ``corners" give rise to
punctures. In computing the Euler characteristic,
 every puncture contributes a -1.

{\proclaim { Definition 1.2 \rm {(Measured train-track)}}}(see [PH]) A {\it { measured
train-track}} $(\tau, \mu)$ consists of a train-track $\tau$, and an
assignment of a non-negative number $\mu (b)$ for each branch $b$ of
$\tau$, so that the following condition holds: For any switch $v$ of
$\tau$,

$ \sum \{ \mu(b) | \ b$ an incoming branch to $v \}= \sum \{ \mu(b) |
\ b$ an outgoing branch to $v \}.$ 

The above condition is called the
{\it{switch condition}}. We also use the term switch condition for a
particular switch $v$.

{\proclaim { Definition 1.3 \rm{($\pi_1$-train-track)} }} (see [BS]) Suppose
$\eurm S=\eurm S^\eurm p_{\eurm g}$ is a surface with $\chi(\eurm
S)=2-2\eurm g-\eurm p<0$. The universal cover of $\eurm S$ then can be
identified with hyperbolic plane $\Bbb H^2$. Fix a polygon $R$ in
$\Bbb {H}^2$ as a fundamental domain for the action of $\pi_1(\eurm
S)$ on $\Bbb {H}^2$. Notice that $R$ is naturally identified with
$\eurm S$ cut open along a number of arcs. Let $\tau$ be a train-track
in $\eurm S$. We call $\tau$ a $\pi_1$-train-track (with respect to
the choice of $R$) if the following conditions hold: If we look at
$\tau$ in the cut-open surface $R$, there is at most one switch point
on each edge of $R$, no switch points in the interior of $R$, and all
the branches are properly embedded in $R$, joining distinct vertices
in $\partial R$. (not necessarily distinct in $\eurm S$.)

{\proclaim {1.4. The Moves}} (see [PH]) We denote by $ {\Cal {MT}}(\eurm S)$ the space of
all measured train-tracks on a surface $\eurm S$, modulo an
equivalence relation which is generated by the following three moves:

$(i)$ Isotopy.

$(ii)$ Right or left {\it{split}} (Figure 1.1).

\midinsert \ForceWidth{5in} \EPSFbox{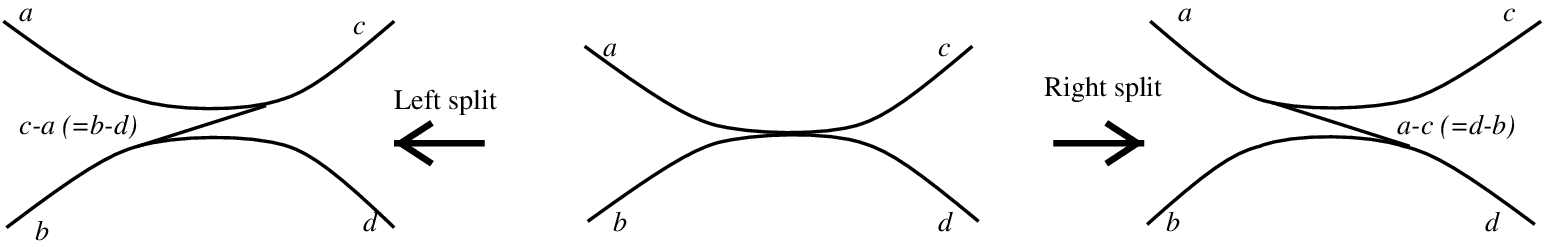} \botcaption{Figure
1.1}{ } \endcaption \endinsert

$(iii)$ {\it{Shift}} (Figure 1.2).

\midinsert \ForceWidth{5in} \EPSFbox{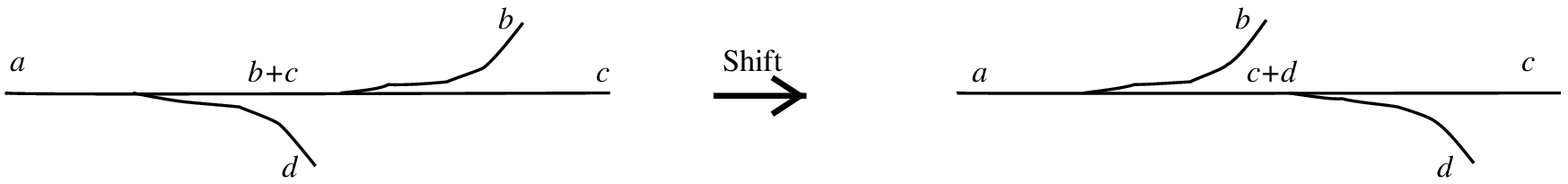} \botcaption{Figure
1.2}{ } \endcaption \endinsert

We have only shown the relevant piece of the train-track in Figures
1.1, 1.2. Notice that the inverse of a split is called a {\it
collapse}.

The set of measures on a train-track $\tau$ is denoted by $V(\tau)$,
and can be identified with a subset of some Euclidean space defined by
a finite set of equalities and inequalities. The set $V(\tau)$ is
closed under (positive) scalar multiplication and addition. In
particular, it is a convex cone.

The following theorem, which is probably due to Thurston,  gives a coordinate system for $\Cal{MT}(\eurm
S)$, in the case which $\eurm S$ has negative Euler characteristic and
is not closed.

\proclaim{ Theorem 1.5} Let $\eurm S$ be a non-closed surface (i.e.,
$\eurm p>0$) with $\chi(\eurm S)<0$, and let $R$ be a polygon
representing a fundamental domain for the action of $\pi_1(\eurm S)$
on the hyperbolic plane. Then any measured train-track on $\eurm S$ is
equivalent to a unique $\pi_1$-train-track with respect to $R$. In particular,
every non-trivial multiple closed curve corresponds to a unique (integral) measured $\pi_1$-train-track.

\endproclaim

This theorem is proved in [HC] (see Theorem 5.1 there) in the case of
a surface with 1 puncture. The general proof is completely
similar.  The following direct corollary gives a piecewise linear
structure on $\Cal{MT}(\eurm S)$.
 
{\proclaim{Corollary 1.6} For a surface $\eurm S$ and polygon $R$ as
above, $\Cal{MT}(\eurm S)$ is the finite union of the cones $V(\tau)$
where $\tau$ ranges over the finite set of $\pi_1$-train-tracks with
respect to $R$.
\endproclaim

For any surface $\eurm S$ the mapping class group $\eurm M_\eurm S$
acts on $\Cal{MT}(\eurm S)$, since if one changes a train-track $\tau$
by any of the moves (i)-(iii) or change a homeomorphism $f:\eurm S \to
\eurm S$ by isotopy, then $f(\tau)$ changes by a sequence of the moves
(i)-(iii).  When
 a homeomorphism $f$ acts on a $\pi_1$-train-track $\tau$ it need not
map it to a $\pi_1$-train-track. Using Theorem 1.5 one can put the
image $f(\tau)$ in the $\pi_1$-train-track by a sequence of the moves
(i)-(iii). We will study how these moves must be performed, and what
the corresponding action of $f$ on $V(\tau)$ is.

Let $\eS=\eS^\ep_\eg$ be a surface with $\chi(\eS)<0$,
 and the polygon 
$R$ be a fundamental domain for the action of $\pi_1(\eS)$ on
$\Bbb H^2$.

{\proclaim{1.7}} Let $n$ be the number of edges in the polygon $R$ and
call the edges $e_1, e_2,..., e_n$ in clockwise order. Give each $e_i$
the orientation induced by the clockwise orientation on $\partial R$.
 If $e_i$ is identified with $e_j$ in $\eurm S$ (obviously with the
 opposite orientation, since $\eurm S$ is orientable), we denote that
by $e_i=e_j^{-1}$.

Pick a base point $x_0$ in the interior of $R$. We want to specify a
set of generators for $\Gamma=\pi_1(\eurm S,
x_0)$. Let $\gamma_i$ be a simple closed
curve based at $x_0$ defined as follows: It starts at $x_0$, it
crosses $e_i$ (it naturally comes out of $e_j=e_i^{-1}$) and then it
goes back to $x_0$, without crossing $\partial R$ any further. The
curve $\gamma_i$ gives rise to an element in $\Gamma$, which by abuse of notation we call $e_i$ too. Notice that the
equation $e_j=e_i^{-1}$ holds in $\Gamma$ as well. It is easy to see
that $e_1,...,e_n$ generate $\Gamma$.

{\proclaim{1.8}} A simple closed curve $C$ can be given by a cyclic
word $e_{\alpha_1}...e_{\alpha_k}$ where $1 \le \alpha_i \le n$. To
draw the curve in $R$ from the given word, just start on the base
point $x_0$, go to $e_{\alpha_1}$, come out of interval
$e_{\alpha_1}^{-1}$ and connect it to $e_{\alpha_2}$, so that it'll
come out of
  $e_{\alpha_2}^{-1}$, etc. All the curves that we consider are
assumed to be tight, i.e., $\alpha_i^{-1} \ne \alpha_{i+1}$ for
 all $i$ (consider $i$ to be a cyclic index modulo $k$).
 
{\proclaim{1.9}} Let's set up some notation for the case when
$\eurm {S=S_g^1}$ is a surface of genus $\eurm g$ with one puncture
$P$, since this case is the simplest case.
 We use the standard fundamental
domain $R$ for the surface $\eurm S$, which is a $4\eurm g$-gon with
edges labeled as $E=E(R)=(a_1,b_1,a_1^{-1},b_1^{-1}...a_\eurm
g,b_\eurm g,a_\eurm g^{-1},b_\eurm g^{-1}),$ in clockwise order. We
call $E$ the edge set.
 
  When we draw curves in $R$, if we are only
interested in their free isotopy class, we draw them off the base
point.  It is important to notice, for example, that the curve given
by the sequence $a_1$ is different from edge $a_1$. It is actually
parallel to the edge $b_1$, but in different orientation. Also, the
curve $b_1$ is parallel to edge $a_1^{-1}$, with the same
orientation. Let's
 introduce the curves $x_1,...,x_\eurm g$. For $1 \le i \le\eurm g$,
the curve $x_i$ is given by the sequence $b_ia_{i+1}$ (take
 the indices mod $\eurm g$, for
example, in the case $i=\eurm g$ in the definition of $x_i$). Let $D_c$
denote the (right-handed) Dehn twist about the simple closed curve
$c$. By [Hu] or [B] we have 
$$\eurm
M_\eurm S=\eM^1_\eurm g=\langle D_{a_1},D_{b_1},...,D_{a_\eurm
g},D_{b_\eurm g},D_{x_1},...,D_{x_{\eurm g-1}}
\rangle.\tag{1.2}$$

One has to notice that, the same set generates $\eM^0_\eg$ if the curves are considered in the closed surface.

 Any mapping class $f$ on $\eS^1_\eg$ ( $\eg >2$)
is specified with its action on the simple
closed curves (with base point)
 $$a_1,b_1,...,a_\eurm g,b_\eurm g.$$
 If so, then for any simple closed curve $c=e_1...e_N$,
$f(c)=f(e_1)...f(e_N)$. This is simply because $f$ induces a
homomorphism on the fundamental group, and if $f$ induces the identity on
$\pi_1(\eS)$, $f$ is the identity mapping class.
(If $\eg=2$ then $f$ also could be hyperelliptic involution.) 

We know that $\eurm M_\eurm S$ is generated by finitely many Dehn
twists. Therefore, it is enough to study the action of a single Dehn
twist on a measured $\pi_1$-train-track $\nu=(\tau,\mu)$.

\topinsert

\ForceWidth{5in} $$\EPSFbox{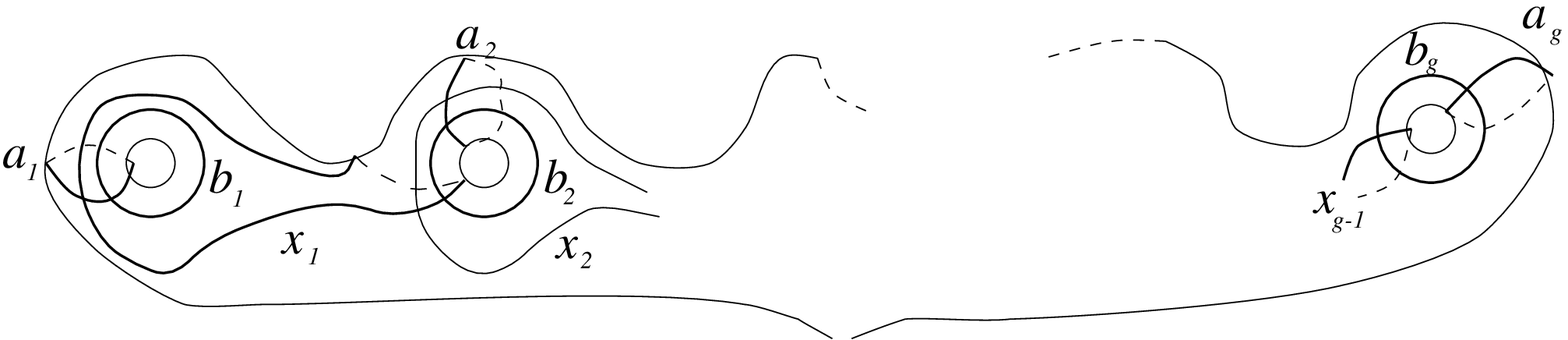}$$ \botcaption{Figure $1.3$}{ }
\endcaption \endinsert

 For a $\pi_1$-train-track $\tau$ on $R$, we call a branch $b$ of
$\tau$ {\it{outer}} if it connects two consecutive edges of the
polygon $R$. Otherwise we call $b$ {\it{inner}}. By out$(\tau)$
(resp. inn$(\tau)$) we mean the set of outer (resp.  inner) branches
of $\tau$. The train-track $\tau$ is identified with the set of
branches of $\tau$.  So $\tau=$inn$(\tau) \ \cup$ out$(\tau)$. We say
a measured train-track $\nu$ is {\it {precisely carried}} on a
$\pi_1$-train-track $\tau$, if $\nu$ is carried on $\tau$ and is not
carried on any sub-train-track of $\tau$.

For a measured $\pi_1$-train-track $\nu=(\tau, \mu)$ ,
 the {\it {total measure}} of $\nu$ is defined by $$T(\nu)=\sum_{b \in
\tau} \mu(b).$$ Notice that $T(a \nu)=aT(\nu)$ for $a>0$. The space of
projective measured train tracks can then be defined by
$$\Cal{PMT}(\eurm S)=\{ \nu \in \Cal{MT}(\eurm S) \ | \ T(\nu)=1 \}.$$
Also the canonical projection $\Cal{MT}(\eurm S) \backslash \{0\}
{\overset{[\cdot ]}\to \longrightarrow} \ \Cal{PMT}(\eurm S)$ can be
defined by $[\nu]=\nu/T(\nu)$.

\

\head $\S 2.$  Complexity of the word problem in
the mapping class groups of  once-punctured surfaces \endhead

Let $\eurm S=\eurm{S^1_g}$. As we saw before,
a generating set for $\eM_\eS$ is given by (1.2).
In this section we consider the following problem: What is the complexity
of computing
(i)  $D_{a_i}(\nu)$ or $D_{b_i}(\nu)$, (ii) $D_{x_i}(\nu)$ and
 (iii) $D_{\tilde a_i}(\nu)$ or $D_{\tilde b_i}(\nu)$ 
for a given 
integral measured $\pi_1$-train-track $\nu=(\tau, \mu)$.
  Let $T(\nu)=\ell$. 
Unfortunately the notation in [HC] is different from our notation.
There $E(R)=(e_1,\cdots,e_{4\eg})$ while here $E(R)=(a_1,b_1,a_1\i,b_1\i,\cdots)$. Also, in [HC], for $1 \le t \le 2\eg$, the curve $b_t$ is defined to be
$e_{2t-1}e_{2t+1}$ for odd $t$ and $e_{2t-2}e_{2t}$ for even $t$. In other words,
our collection of simple closed curves $\{a_1,b_1,\cdots, a_{2 \eg}, b_{2 \eg} \}$ is the same as $\{ b_1,\cdots, b_{2\eg}\}$ in [HC]. To make the notation 
clear, let $b^*_t$ denote the $b_t$ in [HC]. We will only use this notation in 2.1 below. 
Let's look at the complexity of the computation of  $D_{b^*_t}(\nu)$.

{\proclaim{2.1. Complexity of computing $D_{b^*_t}(\nu)$}}

{\bf{1.}} Enter $\nu$ in the machine in the following form:
$L(\nu)=\{(e_i, e_j, \mu(e_i, e_j))\}_{i,j}$, where $e_i,e_j$ are edges of $R$,  
and $\mu_{ij}=\mu(e_i,e_j)>0$
is the corresponding measure. Since there can be at most $2|E(R)|-3$
branches in $\tau$, $L(\nu)$ has $O(\eurm g)$ elements. Since
$1 \le e_i,e_j \le 4\eurm g$ and $1 \le \mu_{ij} \le  \ell$, this has 
complexity $O( \eurm g (\log \ell + \log \eurm g))=O(\eg \log (\eg \ell))$.
Notice that entering a number of size $O(N)$ into the machine has complexity 
$O(\log N)$.

{\bf{2.}} Put $k(i)=2i$ for $i$ odd and $k(i)=2i-1$ for $i$ even.
Check if $\mu_{k(t), k(t+1)}=0$. Looking at $L(\nu)$, this has complexity
$O(\eg \log (\eg \ell))$.

{\bf{3.}}If $\mu_{k(t), k(t+1)} \ne 0$, go to step 5.  If $\mu_{k(t), k(t+1)}=0$,  the resulting train-track after applying 
$D_{b^*_t}$ is collapsible to a $\pi_1$-train-track. One can obtain $D_{b_t}(\nu)$
by changing all $(e_i, e_{k(t)}, \mu(e_i, e_{k(t)}))$ to
 $(e_i, e_{k(t)+1}, \mu(e_i, e_{k(t)}))$ and adding $(e_{k(i)}, e_{k(i)-1},
\sum_i \mu(e_i, e_{k(t)})))$ to $L(\nu)$. This results in a collection 
$L_1=L_1(D_{b^*_t}(\nu))$. Notice that obtaining $L_1$ has complexity 
$O(\eg \log (\eg \ell))$ as well. Also, $|L_1|=O(\eg)$. Also notice that since we added only some of the terms of $L(\nu)$ at most once, $T(D_{b^*_t}(\nu)) \le
2 T(\nu)=2 \ell$.

{\bf{4.}} To obtain $L( D_{b^*_t}(\nu))$ from $L_1$, sort $L_1$ Lexicographically in terms of the first two components. Then combine any string  of consecutive
 terms of the form
$(e, e\p, m_1),\cdots , (e, e\p, m_s)$ to $(e, e\p, \sum_i m_i)$. This gives
$L(D_{b^*_t}(\nu))$, as desired. The sorting and combining processes each
 have complexity $O(\eg \log (\eg \ell))$.

{\bf{5.}} If $\mu_{k(t), k(t+1)} \ne 0$,  the resulting train-track after applying 
$D_{b^*_t}$ is not collapsible to a $\pi_1$-train-track.
 As in step 3, one can obtain $D_{b^*_t}(\nu)$
by changing all $(e_i, e_{k(t)}, \mu(e_i, e_{k(t)}))$ to
 $(e_i, e_{k(t)+1}, \mu(e_i, e_{k(t)}))$ and adding $(e_{k(i)}, e_{k(i)-1},
\sum_i \mu(e_i, e_{k(t)})))$ to $L(\nu)$. This results in a collection 
$L_1=L_1(D_{b^*_t}(\nu))$. Notice that obtaining $L_1$ has complexity 
$O(\eg \log (\eg \ell))$ as well. The list $L_1$ has an element of the form 
$(e_{k(i)+1}, e_{k(i)+1}, \mu(k(i), k(i)+1))$. Drop this from $L_1$. This is 
equivalent to reducing the bad curve. Following 3.2 in [HC],
Now we have to do a split. Create two lists $\{A_1,\cdots ,A_n\}$ and 
$\{B_1,B_2\}$, as instructed in Figures 8 and 9 there. Then decide which split 
to do as in Figure 10. All these steps can be implemented with complexity
$O(\eg \log (\eg \ell))$. Change $L_1$ accordingly, and then go to step 4
to obtain $L(D_{b^*_t}(\nu))$. The estimate  $T(D_{b^*_t}(\nu)) \le
2 T(\nu)=2 \ell$ still holds.

The steps 1-5 show that

\proclaim {Theorem 2.2} Let $\nu=(\tau, \mu)$ be an  integral measured
 $\pi_1$-train track with respect to the standard fundamental domain $R$
 for $\eS^1_\eg$ with $T(\nu)=\ell$. 
Then one can compute $D_{a_t}(\nu)$ and $D_{b_t}(\nu)$ with complexity $O( \eg  \log (\eg \ell))$
and one has $T(D_{a_t}(\nu)) \le 2\ell$ and $T(D_{b_t}(\nu)) \le 2\ell$.
\endproclaim

Similarly, but a  more detailed argument one can obtain from 3.3 in [HC] the following:

\proclaim {Theorem 2.3} Let $\nu=(\tau, \mu)$ be an  integral measured
 $\pi_1$-train track with respect to the standard fundamental domain $R$
 for $\eS^1_\eg$ with $T(\nu)=\ell$. Let $x_t$ be the simple closed curve
$b_ta_{t+1}$.
Then one can compute $D_{x_t}(\nu)$ with complexity $O( \eg  \log (\eg \ell))$
and one has $T(D_{x_t}(\nu)) \le 3\ell$.
\endproclaim

The case of $D_{\tilde a_t}$ and $D_{\tilde b_t}$
was not discussed in [HC]. However, 
similar arguments can be applied. Since 
$T(\tilde a_t)=T(\tilde b_t)=T(\tilde b_t)=4\eg+1$,
 one needs to do steps similar to step 5  in 2.1 
 $O(\eg)$ times, therefore giving:

\proclaim{Theorem 2.4} There are 4 integral measured $\pi_1$-train-tracks $\nu_i$,
$i=1,\cdots,4$ on $\eS^{1}_\eg$, $\eg \ge 2$, such that for $f \in 
 \eM^{1}_\eg$, the following condition implies $f=\id$. 
\roster
\item"(*)"
 $f(\nu_i)=\nu_i$ for  $i=1,\cdots,4$.  
\endroster
\endproclaim

\demo{Proof} Figure $2.1$ shows a ``pair of pant''
decomposition of $\eS^1_\eg$ by a set of simple closed curves
$P=\{\a_i,\b_i, \g_i\}_{i=1}^{\eg-1} \cup \{\d\}$.
For any curve $\rho \in P$ one can define the simple closed curve  $\rho\p$
by Figure $2.2$. If a mapping class $f$ fixes all the curves in $P$, then it must be a product of $D_\rho^{\pm 1}$, $\rho \in P$. If, moreover, $f$ fixes all
$\rho\p$, $\rho \in P$, then $f=$id. Set 
$$\nu_1=\{\a_i,\b_i, \g_i\}_{i=1}^{\eg-1} \cup \{\d\p\},$$
$$\nu_2=\{\g_1\p,\cdots \g_{\eg-1}\p \}  \cup \{\d\},$$
$$\nu_3=\{\a_1\p,\a_3\p, \cdots\} \cup \{\b_2\p,\b_4\p,\cdots\},$$
$$\nu_3=\{\a_2\p,\a_4\p, \cdots\} \cup \{\b_1\p,\b_3\p,\cdots\}.$$
It is easy to see that each collection $\nu_i$ consists of mutually disjoint curves, so can be made into a measured $\pi_1$-train-track. Moreover, by construction, if a mapping class fixes all $\nu_i$, it must be the identity. $\spadesuit$
\enddemo

\midinsert \ForceWidth{5in} 
$$\EPSFbox{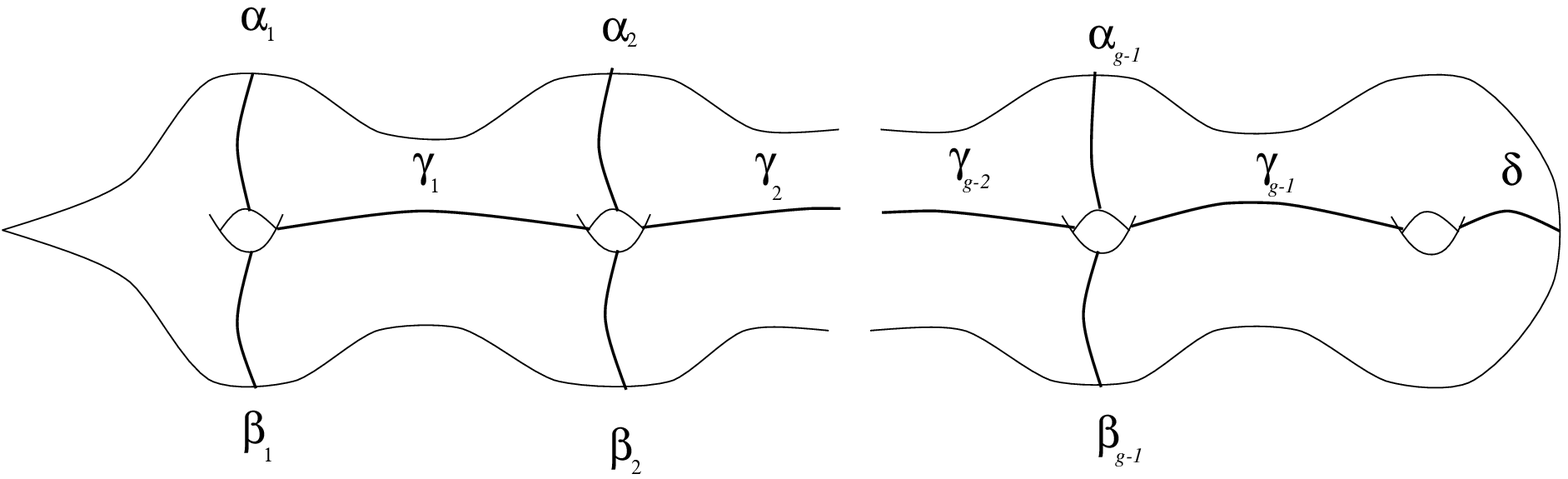}$$ \botcaption{Figure
$2.1$}{ } \endcaption \endinsert

\midinsert \ForceWidth{2in} 
$$\EPSFbox{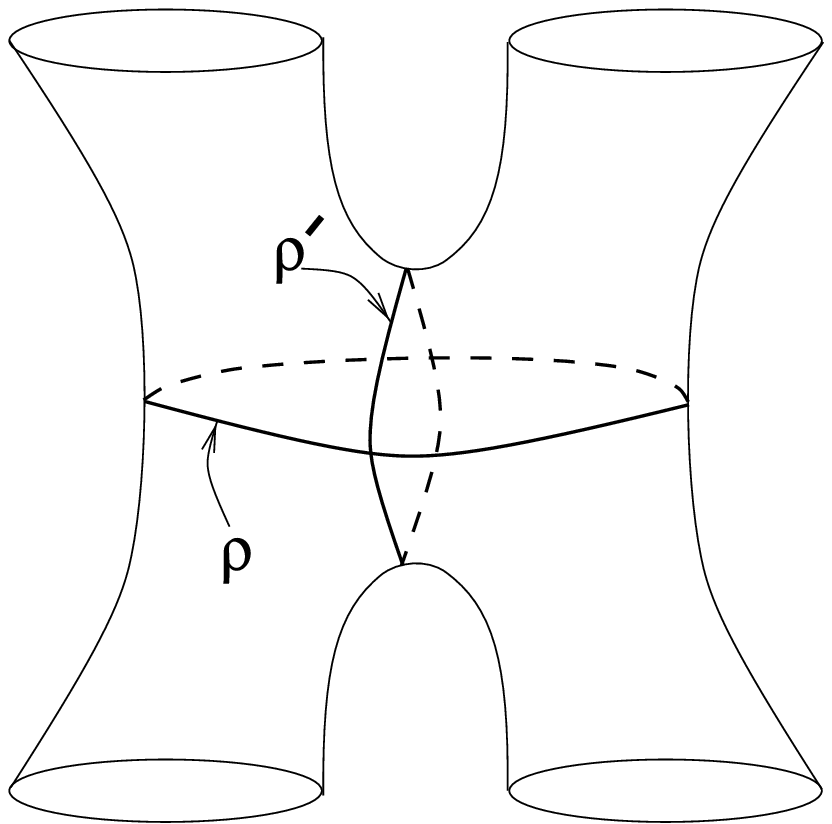}$$ \botcaption{Figure
$2.2$}{ } \endcaption \endinsert

\proclaim{Theorem 2.5} The word problem in $\eM^{1}_\eg$ has complexity 
$O(|w|^2 \eg+|w| \eg \log \eg) $, for a word $w$ in the generators given in
(1.2)
of length $|w|$.

\endproclaim

\demo{Proof}
Let $K=\max\{T(\nu_1),\cdots, T(\nu_4)\}$.
 Notice that $K=O(\eg)$. Compute each $w(\nu_i)$, $i=1,\cdots,4$ by applying 
generators iteratively. At each step, the total measure grows by a factor
of at most $3.$
Therefore the total complexity is 
$$ O( \eg \log(\eg K)+ \eg \log(3\eg K)+ \cdots+ 
\eg \log(3^{|w|-1} \eg K))=
O(|w|^2 \eg+|w| \eg \log \eg).$$
Now check if $w(\nu_i)=\nu_i$. This takes $O(|w|\eg \log \eg)$.
This shows that the word problem in $ \eM^{1}_\eg$ has 
complexity $O(|w|^2 \eg+|w| \eg \log \eg)$. $\spadesuit$
\enddemo

{\proclaim{Conjecture 2.6}} The bound given in Theorem 2.5 is in fact optimal.

\

\head $\S 3.$  The complexity of the word problem in Braid Groups \endhead

To study the complexity of the word problem in the Braid groups $\eB_\en$,
$\en \ge 3$, we can use similar methods as before. First we study the mapping class group
$\eM^{\en+1}_0$ of the $(\en+1)$-punctured sphere $\eS_0^{\en+1}$.
Let's call the punctures $P_0,...,P_\en$. Because of the nature of braid groups, we have
to allow mapping classes to permute the punctures $P_1,...,P_\en$ but keep $P_0$
fixed. Let's call this extended group 
$\tilde \eM^{\en+1}_0$. Then we have an exact sequence
$$1 \to \eM^{\en+1}_0 \to \tilde \eM^{\en+1}_0 \to S_{\en} \to 1,$$
where $S_{\en}$ is the symmetric group on $\en$ elements.

We can use the fundamental polygon  $R=(a_1,a_1\i, \cdots, a_\en, a_\en\i)$
to represent $\eS=\eS_0^{\en+1}$. Let's assume that $P_i$ is the vertex shared
by $a_i,a_i\i$.  We can look at the space of measured train-tracks on $\eS$. As in [HC], one can prove that any measured train-track can be represented 
 uniquely  as a measured $\pi_1$-train-track.

To determine the action of $f \in \tilde \eM=\tilde \eM^{\en+1}_0$ on a measured
$\pi_1$-train-track  $\nu=(\tau,\mu)$ one has to also specify a permutation 
$\sigma \in S_\en$. 
The group $\tilde \eM$ is generated by $\en-1$ half-twists $H_i$  along the curves $\gamma_i=a_ia_{i+1}$ for $i=1,...,\en-1$.

By a half-twist along $\gamma_i$ we mean the following mapping class, which
interchanges $P_i$ and $P_{i+1}$, and is obtained by cutting $\eS$ along a strip
parallel to $\gamma_i$, rotating the component containing $P_i,P_{i+1}$ by 
$180^\circ$, and then gluing to the rest of the surface continuously, twisting 
towards left (we could use twists to right as well, since the situation is completely symmetric). We will use the set of generators $H_1, \cdots H_{\en-1}$
as our basic set of generators for $\eM$.

\midinsert \ForceWidth{3in} 
$$\EPSFbox{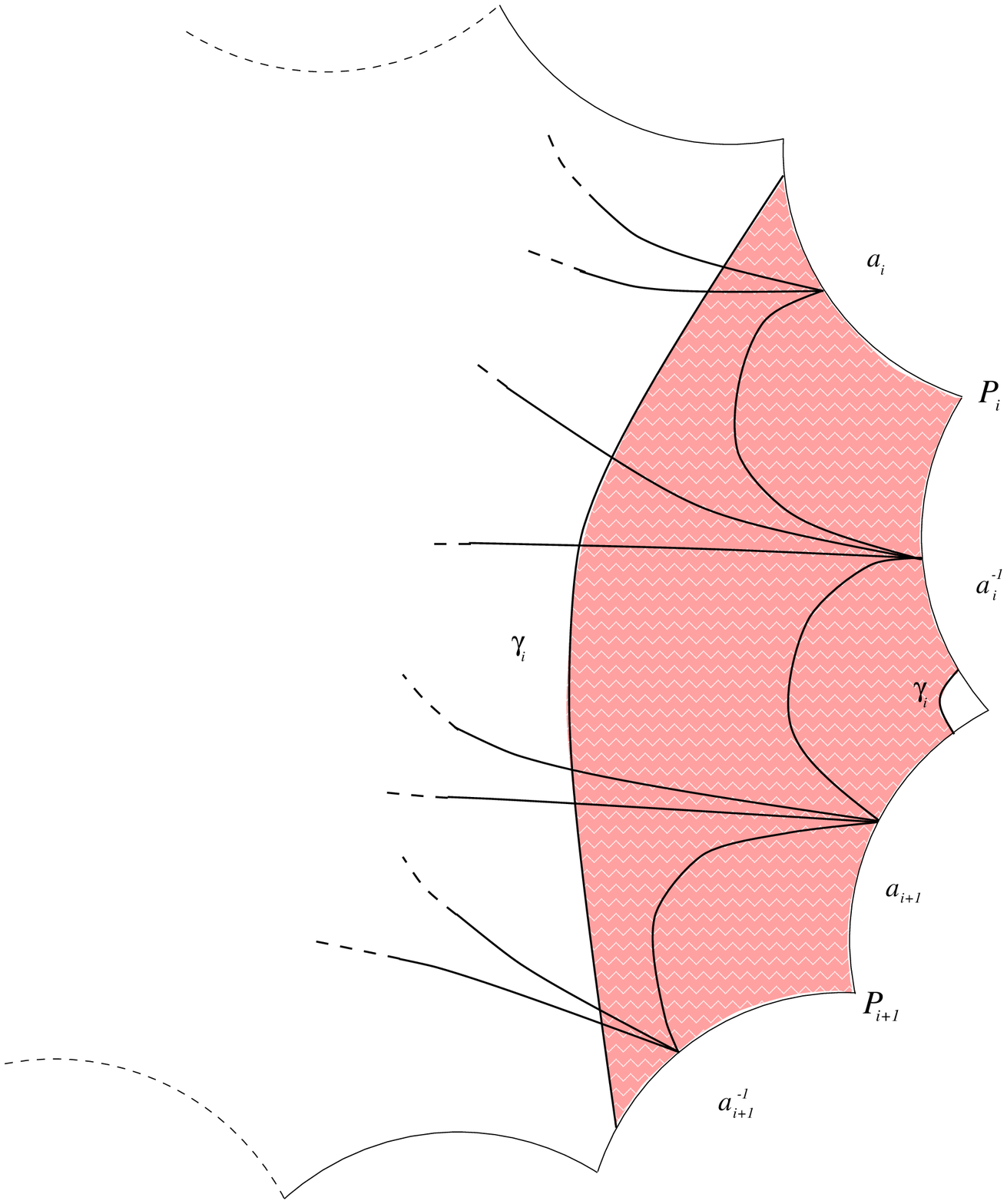}$$ \botcaption{Figure
$3.1$}{ } \endcaption \endinsert

\midinsert \ForceWidth{3in} 
$$\EPSFbox{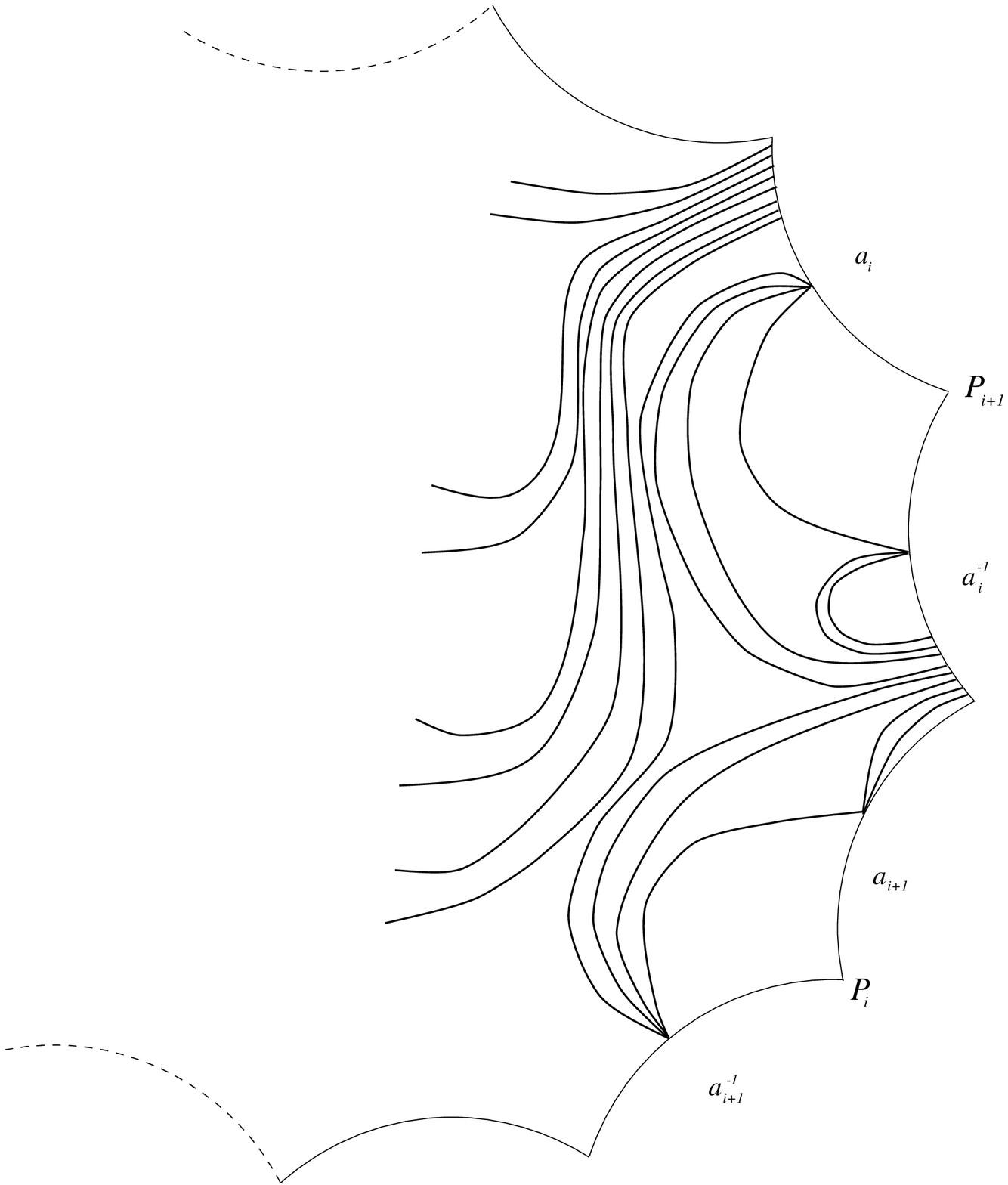}$$ \botcaption{Figure
$3.2$}{ } \endcaption \endinsert

{\proclaim{3.1 Computation of $H_i$ on a measured $\pi_1$-train-track}}

 Now let's see how one can compute
$H_i(\nu)$ for a given measured $\pi_1$-train-track $\nu=( \tau, \mu)$
on $R$. Look at Figure $3.1$, where we have a ``general'' $\pi_1$-train-track.
 We have shaded the region bounded by $\gamma_i$ 
containing $P_i$ and $P_{i+1}$. The outcome of $H_i(\tau)$ is shown in Figure 
$3.2$. To put $H_i(\nu)$ in $\pi_1$-train-track form, we have to consider different cases, as follows:

{\bf{Case 1.}} $\tau$ and $\gamma_i$ do not intersect.
 To get $H_i(\nu)$, we just have to change
 the branches according to the rotation of the 
hexagon bounded by $a_i^{\pm 1}, a_{i+1}^{\pm 1}$
and $\g_i$ by $180^\circ$. Namely, $a_i^{\pm 1} \to a_{i+1}^{\pm 1}$ and 
$a_{i+1}^{\pm 1} \to a_{i}^{\pm 1}$.
This can be done by searching through a list of length $O(\eg)$
 and replacing numbers of order $T(\nu)$.

{\bf{Case 2.}} $\mu(a_{i+1}\i, a_k^{\pm 1}) \ne 0 $
only possibly for $k=i,i+1$. In this case $H_i(\tau)$ is collapsible to a 
$\pi_1$-train-track. Therefore $H(\nu)$ can be computed by $O(\en)$ additions
of numbers $\le T(\nu)$.

{\bf{Case 3.}} Otherwise. In this case there are going to be bad curves, i.e.,
curves going from $a_i\i$ to $a_i\i$. By reducing the bad curves one can see 
that after a split the resulting train-track will be collapsible to 
a $\pi_1$-train-track. Again the number of operations needed to obtain
 the answer is $O(\en)$, and the numbers involved are $O(T(\nu))$.

This finishes the computation. One can observe that this computation is much
less detailed that the corresponding one in $\eM^1_\eg$.  Let's summarize the 
above discussions in the following Theorem:

\proclaim{Theorem 3.2}
Let $\nu=(\tau,\mu)$ be a measured $\pi_1$-train-track on the standard fundamental domain $R$  for $\eS^{\en+1}_0$ with $T(\nu)=\ell$. Let $H_i$ be one of the standard generators of $\tilde \eM^{\en+1}_0$. Then one can compute
$H_i(\nu)$ as a measured $\pi_1$-train-track with complexity 
$O(\en (\log (\en \ell))$. Moreover, $T(H_i(\nu)) \le 2 \ell$.
\endproclaim

The following is similar to Theorem 2.4.

\proclaim{Theorem 3.3} There are 3 integral measured $\pi_1$-train-tracks $\nu_i$,
$i=1,2,3$ on $S^{\en+1}_0$, $\en\ge 3$, such that for $f \in \tilde
 \eM^{\en+1}_0$, the following condition implies $f=\id$. 
\roster
\item"(*)"
 $f(\nu_i)=\nu_i$ for  $i=1,2,3$.  
\endroster
\endproclaim

\midinsert

\ForceWidth{5in}
$$ \EPSFbox{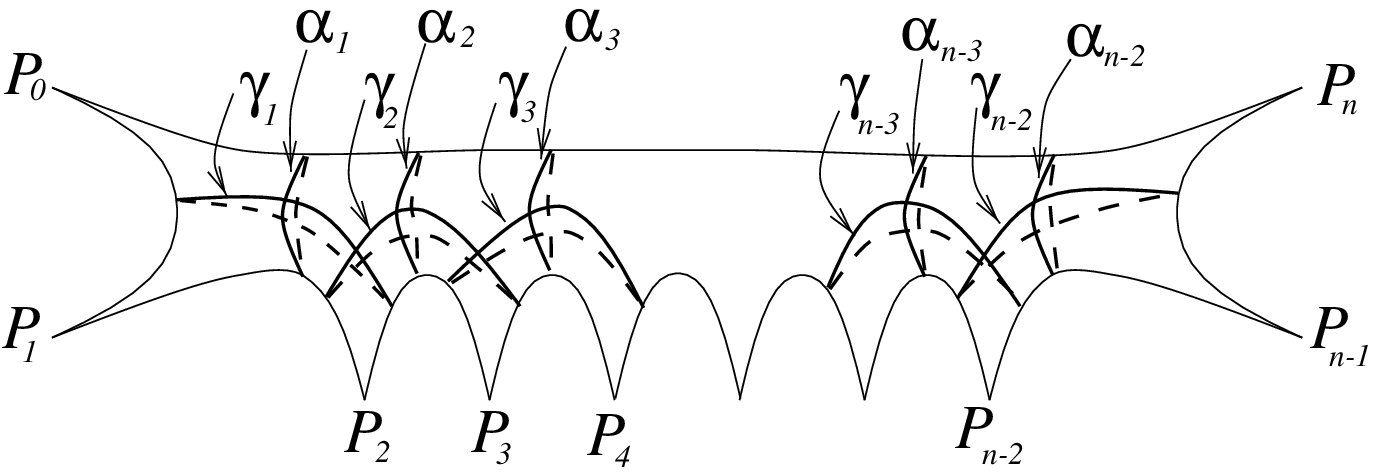}$$ \botcaption{Figure $3.3$}{ }
\endcaption \endinsert

\demo{Proof}
The $(\en+1)$-punctured sphere  can be divided up into ``pairs of pants''
using the simple closed curves $\a_1, \cdots, \a_{\en-2}$. See Figure $3.3$. If $f$ fixes
$\g_1, \cdots , \g_{\en-2}$ and $\a_1, \cdots, \a_{\en-2}$ 
then it has to fix all the punctures. This is easy to see when $\en\ge 4$.
If $\en=3$, i.e., there are 4 punctures, then use the fact that $P_0$
is fixed  by all mapping classes $f \in \tilde \eM^4_0$.
It follows that $f$ must be a 
product of twists in $\a_i$. If $f$ fixes $\g_i$, $i=1,\cdots \en-2$,
then  $f$ can not have a twist in $\a_i$, so  $f=$id.
 Now let $\nu_1$ be the measured train-track
obtained by $\{ \a_1, \cdots, \a_{\en-2} \}$, $\nu_2$ be obtained by
$\{\g_1,\g_3, \cdots \}$ and $\nu_3$ be obtained by $\{\g_2,\g_4, \cdots \}$.
Now if $f$ fixes
 $\nu_1$,$\nu_2$ and $\nu_3$ then it fixes all
$P_i,\g_i, \a_i$. Therefore $f=\id$.  $\spadesuit$
\enddemo

\proclaim{Theorem 3.4} The word problem in $\eM^{\en+1}_0$ has complexity 
$O(\en|w|^2+|w| \en \log \en)$, for a word  $w$ in $\{H_1,\cdots H_{n-1}\}$
of length $|w|$. 
\endproclaim

\demo{Proof}
Let $K=\max\{T(\nu_1),T(\nu_2),T(\nu_3)\}$. 
 Notice that $K=O(\en)$. Compute each $w(\nu_i)$, $i=1,2,3$.
 Each has complexity
$$ O( \en \log(\en K)+ \en \log(\en 2K)+ \cdots+ \en \log(\en 2^{|w|-1}K))=
O(\en|w|^2+|w| \en \log \en).$$
Now check if $w(\nu_i)=\nu_i$. This takes $O(\en \log \en |w|)$.
This shows that the word problem in $\tilde \eM^{\en+1}_0$ has 
complexity $O(\en|w|^2+\en \log \en |w|)$. $\spadesuit$
\enddemo

Now we turn to the word problem in the braid groups. The $\en$-braid group
$\eB_\en$ is given
by the mapping class group of an $\en$-punctured disk, with the possibility
of permuting punctures. Notice that
$$1 \to \Bbb Z \to \eB_\en \to \tilde \eM^{\en+1}_0 \to 1.$$
Also, $\eB_\en$ has the Artin presentation
$$ \align (3.1) \ \ \ \ \ \ \ \ \ \ \ \  \ \ \ \ 
\eB_\en=\langle  \s_1,...,\s_{\en-1}| \
&\s_i \s_j=\s_j \s_i 
, \ \ |i-j|\ge 2, \\
& \s_i \s_j \s_i=\s_j \s_i \s_j , \ \  |i-j|=1 \rangle. \endalign$$
It is easily seen geometrically that $\s_i \to H_i$ in the natural projection
$\eB_n \to \tilde \eM^{\en+1}_0.$ Therefore given a word $w$
in of length $|w|$ one can check if the image of $w$ is the identity in 
$ \tilde \eM^{\en+1}_0$ with complexity 
$O(\en|w|^2+|w| \en \log \en)$.
To solve the word problem in $\eB_\en$, we have to only check the  following:
For a word $w \in$ ker$( \eB_\en \to \tilde \eM^{\en+1}_0 )$, is
$w=\id$? Geometrically, this means that the given word is a twist around the boundary of the disk; i.e., a power of $\Delta$,
where $\Delta$ is the generator of the center of $\eB_\en$. We need to know 
if this power is $0$. For this let's take a look at the fundamental domain $R$.
and the arc $\beta$ connecting $P_i$ to a point in the boundary of
 the disk as in Figure $3.4$. We can find the action of $w$ on $\b$.
 It's natural to 
encode the arc $\b$ as a {\it{measured $\pi_1$-train-track with dead-ends}}.

\midinsert \ForceWidth{5in} 
$$\EPSFbox{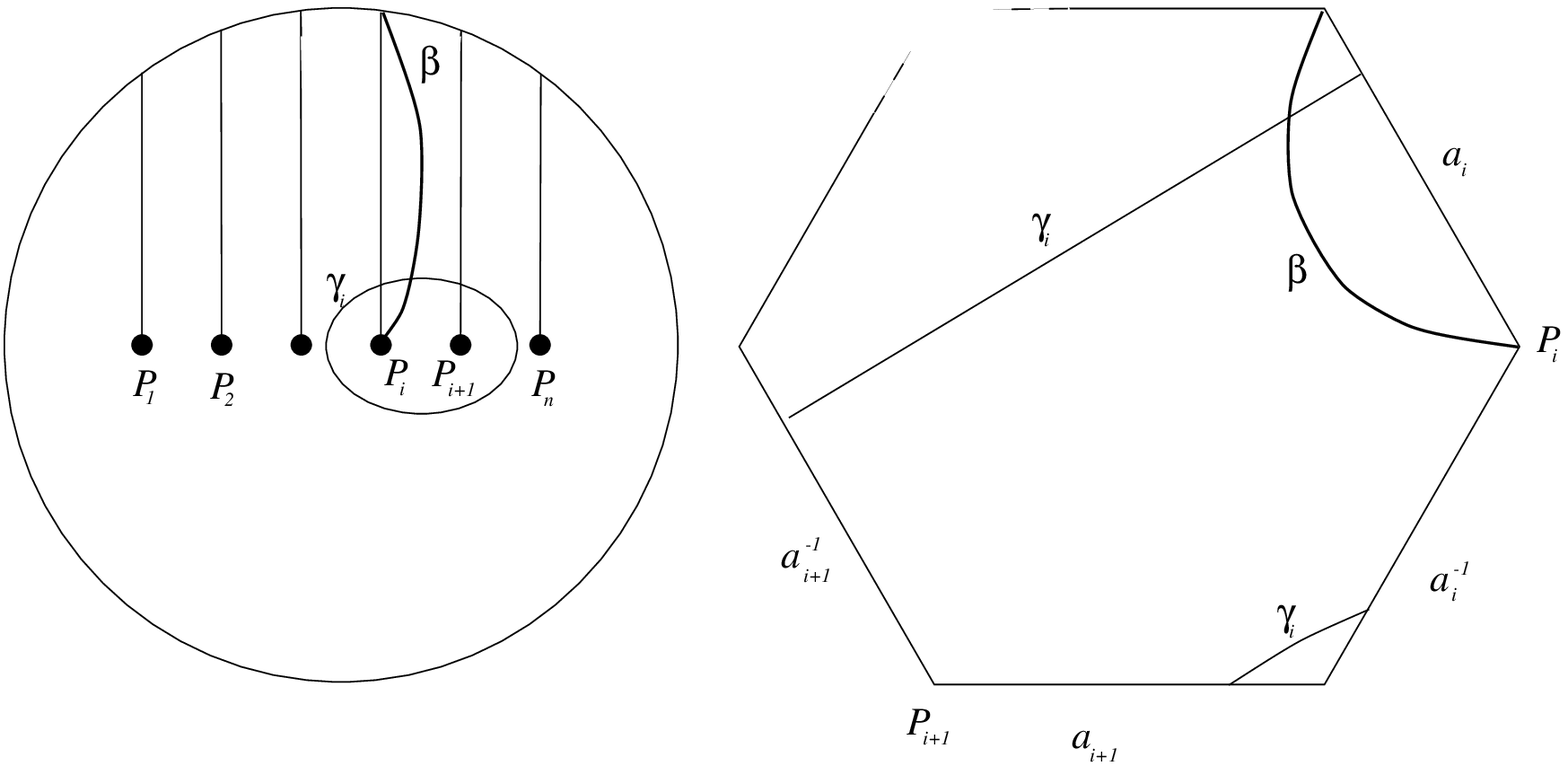}$$ \botcaption{Figure
$3.4$}{ } \endcaption \endinsert

This goes as in the case of measured $\pi_1$-train-tracks, but one has to keep
a neighborhood of the both ends of $\beta$ fixed while applying each generator.
The details are similar to the case of $\pi_1$-train-tracks. In particular,
 this has complexity $O(\en|w|^2+ \en \log \en |w|)$. This implies

\proclaim{Theorem 3.5}
The word problem in the Braid group $\eB_\en$ has complexity
 $O(\en|w|^2+|w| \en \log \en)$, where $|w|$ is the length of the word $w$ in
the Artin generators $\{ \s_1, \cdots, \s_{\en-1} \}$.
\endproclaim

\proclaim{Corollary 3.6}
If $w$ is a word in the Artin generators of  $\eB_\en$  of length $|w|$,
with $|w| \ge \log \en$, one can determine if $w=\id$  with complexity
$O(|w|^2 \en)$ {\it{on a Turing Machine}}.
\endproclaim

Using similar ideas with a standard fundamental domain for the surface 
$\eS^\ep_\eg$ with $\ep\ge 2$ incorporating the cases of once-punctured surfaces and braid groups one can similarly prove:

\proclaim{Theorem 3.7} The word problem in $\eM^{\ep}_\eg$ has complexity 
$O(|w|^2(\eg+\ep)+|w|(\eg+\ep)\log (\eg+\ep))$, for a word $w$ in a set of
``standard''  generators
of length $|w|$.
\endproclaim

\

\head $\S 4.$  The case of a closed surface \endhead

 Now let's discuss the case of a closed surface, i.e., when $\eurm
S=\eurm S^0_\eurm g= \eurm S_\eurm g$. The basic group structure of
$\eurm{M_g}$ in terms of $\eurm{M^1_g}$ is given by the short exact
sequence (see [B]) $$1 \to \pi_1(\eurm S, *) \to \eurm{M^1_g} {
\overset{\phi}\to{\longrightarrow}} \ \eurm{M_g} \to 1.$$ Here the
canonical map $\phi$ is defined by just forgetting the puncture.

Notice that by this exact sequence the generators in (1.2) can be naturally considered as generators of $\eM^0_\eg$.

Let's introduce an artificial puncture on $\eurm S$; i.e., let's fix a
point $P$ on $\eurm S$, and call the corresponding once-punctured
surface $\eurm S^P$. Also let $R$ be the standard fundamental domain
for $\eurm S$, having all vertices equivalent to $P$ on $\eurm
S$. Suppose a simple closed curve $\alpha$ is given on the closed
surface $\eurm S$, and $\alpha$ does not pass through $P$. The curve
$\alpha$ can be considered as a curve on the punctured surface $\eurm
S^P$, and can be given by a cyclic word $w=e_{s_1}...e_{s_n}$ where
$e_i \in E(R)$, as in 1.8. Notice that by isotoping $\alpha$ in
$\eurm S$, we might obtain a shorter cyclic word. We want to discuss
here a geometric analog of Dehn's well-known algorithm (see [J],
e.g.)  to get
 a shortest representative for $\alpha$.  The shortest representative
is not unique, as we will see below.

Represent $w$ by a measured $\pi_1$-train-track $\nu$ carried
precisely on a $\pi_1$-train-track $\tau$. Let $m=|E(R)|/2=2\eurm
g$. If there is no path $\eurm b=(b_1,...,b_k)$ of outer branches in
$\tau$ such that $k\ge m$ then we claim that $w$ is a shortest
representative for $\alpha$.  Recall that in this case $\alpha$ is
given by a word $w=e_{s_1}...e_{s_n}$ in letters in $E(R)$
representing a simple closed curve $\alpha_0$ on $\eurm S^P$. We can
assume $w$ does not have back-tracking; i.e, $e_{s_i}
 \ne e_{s_{i+1}}\i$ for all $i$ (mod $n$).  If $w$ is not a
 shortest representative for $\alpha$, then there is another word
$w\p=e_{s\p_1}...e_{s\p_{n\p}}$ with $n\p <n$ representing a curve
$\alpha_1$ in $\eurm S^P$ which is isotopic to $\alpha$ in $\eurm
S$. Take an isotopy $\alpha_t$ between $\alpha_0$ and $\alpha_1$ on
$\eurm S$.  By changing $\alpha_t$ a little bit, one can subdivide
this isotopy to subintervals in which, either (i) no part of $\alpha$
is in a small neighborhood of $P$, or (ii) only one segment of
$\alpha$ is passing through $P$, and everything else is fixed. Notice
that in intervals of type (i) the word representing the curve in
$\eurm S^P$ does not change since that part of the isotopy can be
looked at as an isotopy of $\eurm S \backslash \{P\}$
 and $\pi_1(\eurm S \backslash \{P\})$ is
free. Therefore, one can find a finite sequence
$\alpha_0=\alpha_{t_0},\alpha_{t_1},...,\alpha_{t_\ell}=\alpha_1$
which give all the different simple closed curves that appear on
$\eurm S^P$. Every element in this sequence is obtained by the
previous one by taking a piece of $\alpha$ and passing it through
$P$. We can assume that $\alpha_{t_i} \ne \alpha_{t_j}$ for $ i \ne
j$, otherwise we can just drop the repeating part of the
isotopy. Since by assumption $\alpha_0$ does not have a path of outer
branches $\eurm b=(b_1,...,b_k)$ with $k \ge m$, we must have
$T(\alpha_{t_1}) \ge T(\alpha_{t_0}),$ with equality only in the case
in which $\alpha_{t_0}$ has a path $\eurm b$ as above with $k=m-1$,
and the move is to just push the path to the other side of $P$ (Figure
4.1). Notice that Figure 4.1 is drawn in the universal cover of
$\eurm S$.

 If the sequence
$\alpha_{t_0},\alpha_{t_1},...,\alpha_{t_\ell}$ only consists of moves
which push a path of length $m-1$ across the puncture,
 then $T(\alpha_0)=T(\alpha_1)$; i.e., $n=n\p$, which is a
contradiction.  Otherwise, let $1 \le j \le \ell$ be such that $T(
\alpha_{t_0})=...=T(\alpha_{t_{j-1}})<T(\alpha_{t_j})$.  If $j=\ell$
then $n\p>n$ which is a contradiction. Since $\alpha_{t_{j+1}} $ can
not be equal to any of the preceding $\alpha_{t_i}$, it is easy to see
that $T$ keeps monotonically increasing on the sequence
 $\alpha_{t_0},\alpha_{t_1},...,\alpha_{t_\ell}$.
 This shows that $n\p>n$, which is again a contradiction.

 \topinsert \ForceWidth{5in} \EPSFbox{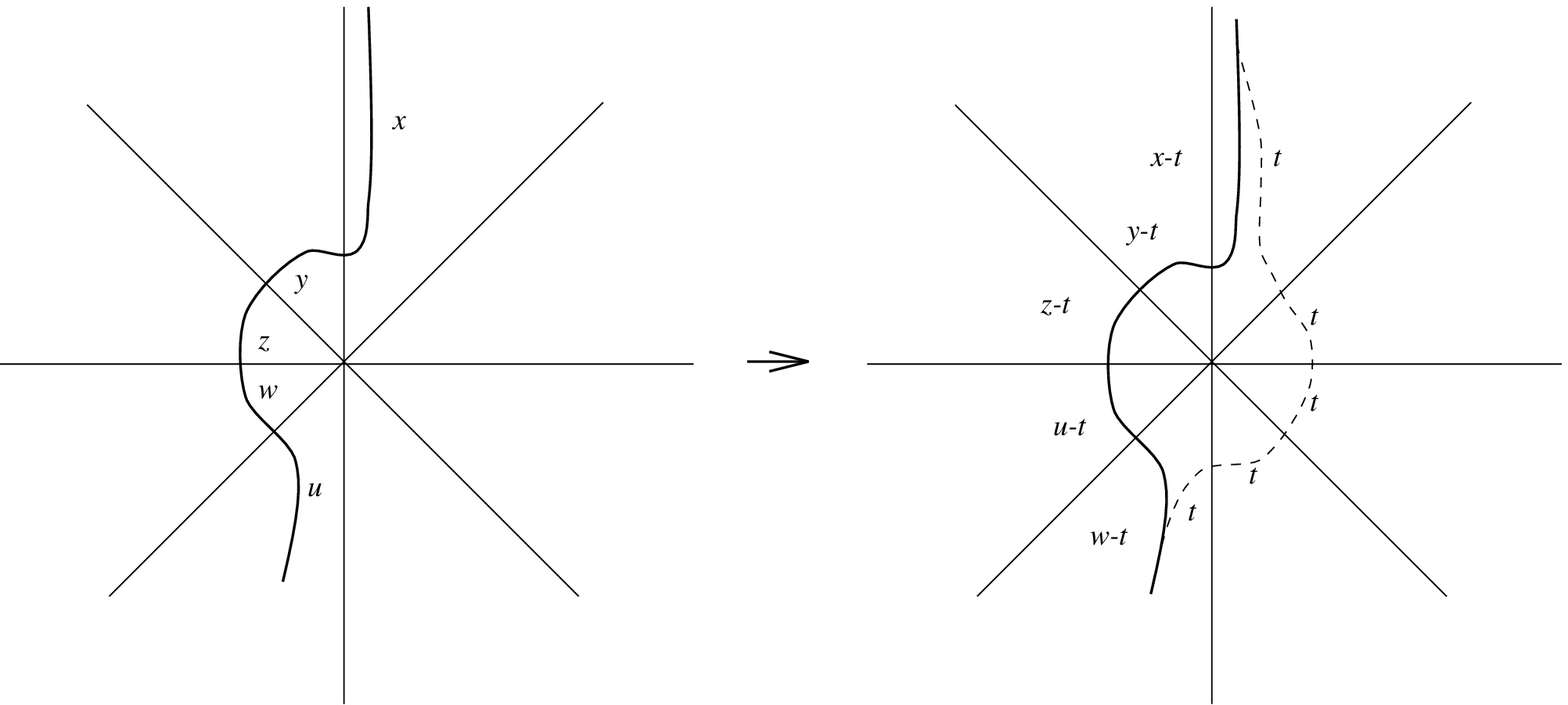}
\botcaption{Figure 4.1}{} \endcaption \endinsert

Let's summarize the above arguments in the following Theorem. A
subword of a word $e_1 \cdots e_n$
 is any word of the form $e_i e_{i+1}\cdots e_j$.

\proclaim{Theorem 4.1} Let $R$ be a standard fundamental domain for
the closed surface $\eurm S=\eurm {S_g}$ with $2m$ edges ($m=2 \eg$),
 with the
vertices of $R$ equivalent to a point $P$ on $\eurm S$. Let $\eurm
S^P$ be a once-punctured surface obtained by fixing $P$ on $\eurm
S$. Let $\alpha$ be a simple closed curve on $\eurm S$ not passing
through $P$, and let $w=e_{s_1}...e_{s_n}$ be a cyclic word in letters
in $E(R)$ representing $\alpha$ up to isotopy in $\eurm S$. Then $w$
is a shortest representative if and only if \roster \item $e_{s_i} \ne
e_{s_{i+1}}\i$ for all $i$ mod $n$, and \item $w$ does not have a
subword of length $\ge m$ consisting of outer branches.  \endroster
Moreover, any two shortest length representatives of $\alpha$ are
related to each other by pushing a finite number of identical subwords
of length $m$ of the outer branches to the other side of
$P$. $\spadesuit$ \endproclaim

With the same assumptions on the fundamental domain $R$, let $\nu$ be
a measured train-track carried precisely on a $\pi_1$-train-track
$\tau$ on $\eurm S$. As we know by now from simple closed curves, the
$\pi_1$-train-track representative is not unique in $\eurm S$. We want
to describe an algorithm to put $\nu$ in a $\pi_1$-train-track form
which has the smallest $T$. We call $\tau$ a
{\it{ reduced-length $\pi_1$-train-track}} if it has no path of outer branches
of length $\ge m$.

\proclaim{Lemma 4.2} If $\nu=(\tau, \mu)$ is a measured
$\pi_1$-train-track on $\eurm S$, there exists a measured
$\pi_1$-train-track $\nu\p=(\tau\p, \mu\p)$ which represents $\nu$ and
it has the smallest possible $T$.  \endproclaim

{\demo {Proof}} Let $\{c_n\}$ be a sequence of simple closed curves on
$\eurm S$ and $\lambda_n>0$ be such that $\lambda_nc_n
 \to \nu$ as $n \to \infty$.  Put each $c_n$ in a reduced form $\tilde
c_n.$ By passing to a subsequence we can assume all the $\tilde c_n$
are carried on a reduced-length $\pi_1$-train-track $\tau\p$. Now one
can look at the sequence $\{ [\tilde c_n] \}$ in $\Cal{PMT}(\eurm
S^P)$. By compactness, this sequence has a convergent
subsequence. Without loss of generality, let's assume $[\tilde c_n]
\to \nu\p \in \Cal{PMT}(\eurm S^P)$. Notice that $\nu\p$ is of reduced
length since it is carried on $\tau\p$.  Using the surjection
$\Cal{PMT}(\eurm S^P) \to \Cal{PMT}(\eurm S)$, one gets a
corresponding convergent sequence $[c_n] \to \nu\p$ in
$\Cal{PMT}(\eurm S)$. We denote the limit point with the same notation
since it is given by the same measured $\pi_1$-train-track. This shows
that $[\nu]= \nu\p$ i.e., $\nu$ is equivalent to a reduced-length
measured $\pi_1$-train-track. Now we have to prove that $T(\nu\p)$ is
minimal among all $T(\nu\pp)$, where $\nu\pp$ is a measured
$\pi_1$-train-track representative for $\nu$. Suppose
$T(\nu\pp)<T(\nu\p)$, for such a $\nu\pp$. Then by definition of the
space of measured train-tracks, there is a finite sequence
$$\nu\p=\nu_1 \to \nu_2 \to...\to \nu_k \to ... \to
\nu_n=\nu\pp\tag{4.1}$$ where each $\nu_j$ is obtained by performing
one of the following moves on $\nu_{j-1}$: (i) Split, (ii) Shift,
(iii) Isotopy without crossing $P$, (iv) Pulling a branch from one
side to the other side of $P$, and (v) Collapse.  It is easily seen
that one can arrange the sequence (4.1) so that $\nu_1,...,\nu_k$ are
obtained by performing the moves of type (i)-(iv), and the rest of
$\nu_j$ are obtained only using the collapse move. Choose a simple
closed curve $c\p$ and $\lambda>0$ such that $c\p$ is carried on
$\tau\p$, and it stays $\e$-close to $\nu_i$ at each step along the
sequence $\nu_1 \to ... \to \nu_k$, as we perform the corresponding
move on $\lambda c\p$,
 where $\e>0$ is an arbitrary pre-chosen number. Here $\e$-close is
used in the sense that at each stage, the sum of the differences the
measures in corresponding branches is bounded above by $\e$. In
particular, $|T(\nu\p)- T(\lambda c\p)|<\e$. After collapsing to
$\nu\pp$, we get a (measured) simple closed curve $\lambda c\pp$ which
is $\e$-close to $\nu\pp$. In particular,
 $|T(\nu\pp)- T(\lambda c\pp)| <\e$. If we choose $2\e<
T(\nu\pp)-T(\nu\p)$, we get $T(\lambda c\pp) < T(\lambda c\p)$, which
contradicts Theorem 4.1, since $c\p$ is carried on a reduced-length
measured $\pi_1$-train-track.  This finishes the proof of the
lemma. $\spadesuit$.

\proclaim{Corollary 4.3} If $\nu=(\tau, \mu)$ is a measured
$\pi_1$-train-track on the closed surface $\eurm S$ carried on a
reduced-length measured $\pi_1$-train-track $\tau$, then $T(\nu)$ is
minimal among $T$ of all other measured $\pi_1$-train-track
representatives of $\nu$.  \endproclaim

A similar limit argument as in the proof of the lemma shows that:

\proclaim{Corollary 4.4} Any two reduced-length measured $\pi_1$
train-track representatives of the same measured train-track on the
closed surface $\eurm S$ are related by the following move: Pulling
some measure off a path of outer branches of length $m-1$, where
$|E(R)|=2m$, to the other side of the puncture $P$.  \endproclaim

Here is an algorithm to put a given measured $\pi_1$-train-track in
the reduced (shortest) form.
 Let's start with a measured $\pi_1$-train-track
$\nu_1=(\tau_1,\mu_1)$ which is not reduced-length. So there is a
unique maximal path $\eurm b=(b_1,..., b_k)$ of outer branches in
$\tau_1$ where $k \ge m$.
  Let $x_i=\mu(b_i)$ be the measure on each branch $b_i$,
$i=1,...,k$. We have to use a move as illustrated in Figure 4.2 to
put $\tau_1$ in a position with smaller $T$. To be able to do that
move, we have to assume $x_i=$min$\{ x_1,..., x_k\}$.  \topinsert
\ForceWidth{4.25in} $$\EPSFbox{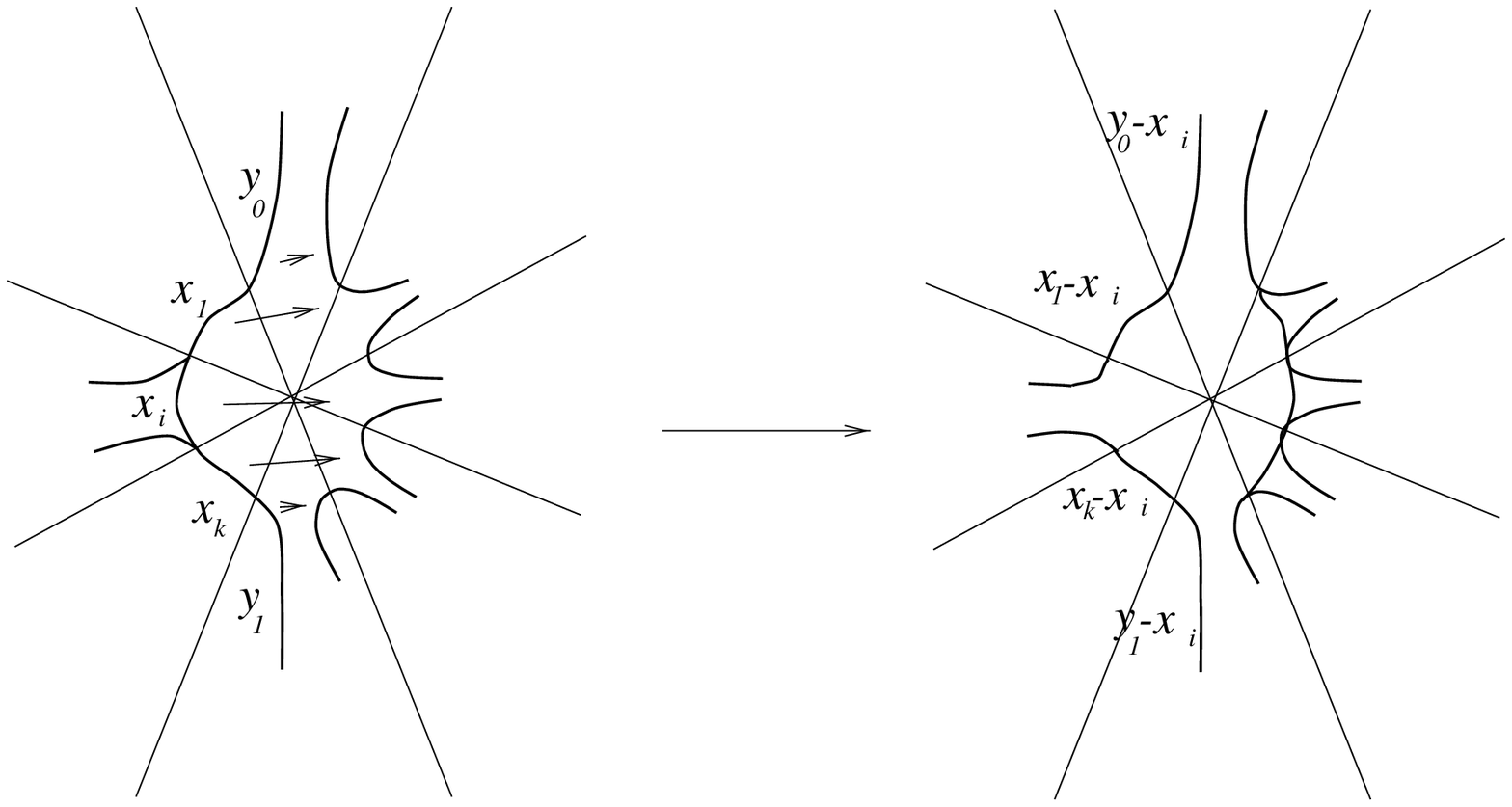}$$ \botcaption{Figure
4.2}{} \endcaption \endinsert

We claim that, after doing the move finitely many times, we will get a
sequence of measured $\pi_1$-train-tracks $\nu_1,...,\nu_t$ where
$\nu_i=(\tau_i, \mu_i)$, and $\tau_t$ is of reduced-length. The reason
is that first of all we know that there is a sequence of moves of type
(i)-(v) putting $\nu_1$ is reduced form.  Now notice that as in the
case of simple closed curves,
 If you make a move and increase $T$, to reduce $T$ later on you have
to undo the move. This proves that there is a sequence to
monotonically decrease $T$, which proves our assertion, since at any
given stage, there is only one way to reduce the $T$, if the
train-track is not already in the reduced-length position.

The analog of Theorem 1.5 is

\proclaim{Theorem 4.5} Let $\eurm S=\eurm S_\eurm g$ where $g\ge 2$,
and let $R$ be an standard fundamental domain for the action of
$\pi_1(\eurm S, *)$ on $\Bbb H^2$. Then every measured train-track is
equivalent to some measured $\pi_1$-train-track $\nu=(\tau,\mu)$ with
respect to $R$ having the smallest possible $T$.  This representative
is unique if and only if $\tau$ has no path of outer branches of
length $|E(R)|/2-1$. Otherwise any representative is obtained from any
other representative by pulling some measure from a path of outer
branches of length $|E(R)|/2-1$ to the other side of the
puncture. $\spadesuit$ \endproclaim

\

\head $\S 5.$ The complexity of the word problem in the mapping class
groups of closed surfaces \endhead

Since the $\pi_1$-train-track representation is not unique for closed surfaces,
the main issue here is the following problem:

{\proclaim{5.1 Problem}} Find the complexity of the following computation: Given an integral measured
$\pi_1$-train-track $\nu=(\tau,\mu)$ on the standard fundamental domain $R$
for the surface $\eS=\eS_\eg$ with $T(\nu)=\ell$, 
compute a $\nu\p=(\tau\p,\mu\p)$ of reduced form such that $\nu\p$ is equivalent to $\nu$ on $\eS$.

Recall that $m=2\eg=|E(R)|/2$. It is easy to check if $\nu$ is not of reduced
 length with complexity $O(\eg)$. One has to check if there is a path of outer
 branches of length $\ge m$. Therefore suppose $\nu$ is not of reduced length,
 to start with.  Let $\eb=(b_1,\cdots, b_{n(\nu)})$ be the unique
 maximal path of outer branches in $\tau$ of length $n(\nu) \ge m$, and let 
$\psi(\nu)=\min\{ \mu(b_1),
\cdots, \mu(b_n)\}$.

We use $n(\nu)$ as a measure of complexity. Notice that
$n(\nu) \le |\out(\tau)|$ 
(Recall that $\out(\tau)$ is
 the set of outer branches of $\tau$). We will put $\nu$ in the 
reduced-length form by a sequence of moves each of which reduces the
 complexity function
$n(.)$. Notice that $\nu$ is of reduced form if $n(\nu)\le 2\eg-1$.
Moreover, it is always possible to reduce $\nu$ such that $|\out(\tau)|
 \le 4\eg-3$,
as we will see below.

{\bf{Case 1. $n(\nu)<4\eg-1$.}}
 Let $x=\mu(b_i)=\psi(\nu)$. We can pull a measure of $x$ to 
the other side of the puncture. This may involve changing some inner branches
 which connect to the both ends of the path $\eb$ to outer ones. 
In particular, this may add a measure of $x$ to at most two of the
 branches in $\eb$.
If none of these branches are $b_i$, then we have reduced the complexity
function, because one can easily see that the added outer branches can not 
extend $\eb$ from either side. Now let's consider the case which pulling
 the measure adds to $b_i$, so that after the pulling, we still have
 $\mu(b_i)=x$. This subtracts $x$ from all the branches of $\eb$ except 
for $b_i$ and possibly another branch $b_j$. By examining the size of the 
measures $\mu(b_k)$, $1\le k \le n(\nu)$ and the ones connecting 
to the endpoints
of $\eb$, we can see how many times this move is possible, and we can do them
all at once. After we do that, there is a $k\ne i,j$ such that
 $\mu(\b_k)$ has become $<x$, which means $\psi(\nu)$ is now $<x$. Now pull this measure across the puncture, and this will reduce the complexity function.
This shows that one can put $\nu$ in reduced-length form after $O(\eg)$ steps.
Each step involves $O(\eg)$ operations on numbers which are $O(T(\nu))$.
Therefore, the complexity of putting $\nu$ in reduced-form in  this case is 
$O(\eg^2 \log T(\nu))$. If at the end the final $\nu$ satisfies
 $n(\nu)=
2\eg-1$, then one can easily force $|\out(\tau)|\le 4\eg-3$:
If  $n(\nu)=
2\eg-1$ and $|\out(\tau)|=4\eg-2$, by finding the outer
 branch with smallest measure and pulling that measure through the puncture in a similar fashion as above, we get $|\out(\tau)| \le 4\eg-3$. 

{\bf{Case 2. $n(\nu)=4\eg-1$.}} (Equivalently,
$|\out(\tau)|=4\eg-1$.) In this case the complexity of the problem 
can be much higher, in fact it will be of linear order with respect to 
$T(\nu)$. The problem is that one can pull a small piece of the curve $\nu$ 
around
arbitrarily long and then hook it up with the puncture.
 Then to simplify the curve one has to undo that, 
which has complexity $O(T(\nu))$. See Figure $5.1$.

\midinsert \ForceWidth{4in} 
$$\EPSFbox{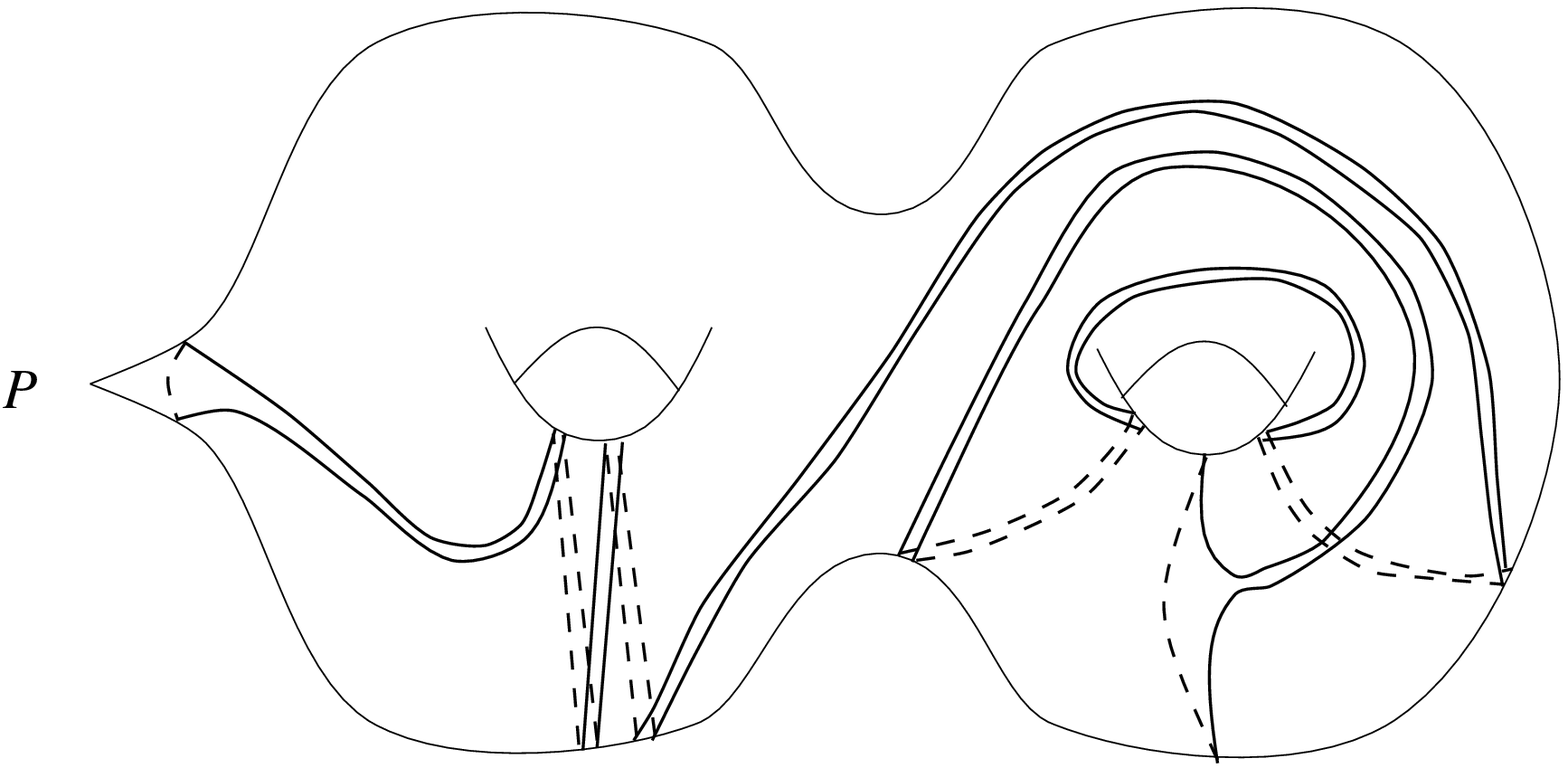}$$ \botcaption{Figure
$5.1$}{ } \endcaption \endinsert

{\proclaim{5.2. Solution to the word problem}}
In the solution to the word problem in $\eM^0_\eg$ we have to avoid Case 2
in 5.1, because it will have an effect of making it exponential, since our polynomial algorithms are all based on the fact that the computations with a curve
are of
order $\log N$, if the size of the curve at hand is $N$.

Here is our strategy for the solution of the word problem in $\eM^0_\eg$:
Let $w=h_1 \cdots h_n$ be a word in the basic set of generators of
$\eM^0_\eg$ (see (1.2) and the note below it). similar to Theorem 2.4, We know that there are 4 measured $\pi_1$-train-tracks
$\nu_1, \cdots, \nu_4$ with $T(\nu_i)=O(\eg)$ on $\eS_\eg$ such that if
$w(\nu_i)=\nu_i$ for all $1 \le i \le
4$
then $w=$id. (This holds only for $\eg \ge 3$; in $\eM_2$ there is 
a mapping class of order 2 fixing all simple closed curves). Put $\nu^{(0)}_i=\nu_i$ 
and $\nu^{(j+1)}_i=h_{n-j}(\nu^{(j)}_i)$. Notice that $\nu^{(n)}_i=w(\nu_i)$.
For $j=0,
\cdots, n$, we compute $\nu^{(j)}_i$. After each computation, we put
$\nu^{(j)}_i$ in the reduced-length form. What we would like to show is that, 
if $h$ is a generator and $\nu$ is of reduced length, $h(\nu)$ can be put
into reduced-length form with complexity $O(\log T(\nu))$ with respect to
$T(\nu)$. For that we have to again look closely how each of the generators
act on a reduced-length measured-train-track $\nu=(\tau, \mu)$. By the above argument
in Case 1, it is enough to show that $n(h(\nu))<4\eg-1$.

\proclaim{Lemma 5.3}
Suppose $h^{\pm1}$ is one of the generators in (1.2), and $\nu=(\tau,\mu)$
 is an
integral measured $\pi_1$-train-track on the standard fundamental domain for
$\eS_\eg$, $\eg \ge 2$ of reduced-length. Put $h(\nu)=(\tau_1,\mu_1)$.
Then $n(\tau_1) <4\eg-1$, or equivalently $|\out(\tau_1)| <4\eg-1$  .
\endproclaim

{\demo{Proof}}
We will only discuss the cases which $h$ is a generator in (1.2).
The cases where $h\i$ is a generator are done by symmetry.

{\bf{Case 1.}} $h=D_{a_t}$.
Let $\nu=(\tau,\mu)$ be a reduced-length measured $\pi_1$-train-track
with $n(\nu) \le 2\eg-1$ and
$|\out (\tau)|\le 4\eg-3$  (see the argument in Case 1 in 5.1 above).
We claim that $n(\nu_1)<4\eg-1$. The proof has many steps.

(i) $\mu(a_t\i,b_t)=0$. (No bad curves) Notice that  $\mu_1(a_t\i,b_t)=0.$
If $\eg>2$, at least one of $\mu(a_{t+1}, b_{t+1}),
\cdots, \mu(a_{t-1}\i, b_{t-1}\i)$ must be 0, which stays 0 with 
 $\mu_1$ instead
of $\mu$. This shows that $n(\nu_1)<4\eg-1$. Suppose  $\eg=2$ and, say $t=1$.
 Since
$\nu$ is of reduced-length, one of the values
$$\mu(b_1\i,a_2),\mu(a_2,b_2),\mu(b_2,a_2\i),\mu(a_2\i,b_2\i)$$
must be 0 and stays 0 if we replace $\mu$ by $\mu_1$. Therefore the estimate
$n(\nu_1)<4\eg-1$ holds in this case too.

(ii) $\mu(a_t\i,b_t)\ne 0$ but $\mu(\a_t,b_t)=0$.
Then again $\mu_1(a_t\i,b_t)=0$ and the argument is similar to (i).

(iii) $\mu(a_t\i,b_t)\ne 0$ and  $\mu(\a_t,b_t)\ne 0$.
In this case $\out(\tau)=\out(\tau_1)$ unless
$\mu(b_t, b_t\i) \ne 0$ and $\mu(a_t\i, b_t\i)=0$, in which case 
$\out(\tau_1)=\out(\tau) \cup \{(a_t\i,b_t\i)\}$. Since $n(\tau)<4\eg-2$,
$n(\tau_1)<4\eg-1$.

{\bf{Case 2.}} $h=D_{b_t}$.
This case is similar to case 1.

{\bf{Case 3.}} $h=D_{x_t}$.

(i) No bad curves. This means that
$\mu(b_t,e)=0$ for $e \in E(R)\backslash  \{a_t\i,b_t\i, a_{t+1},
b_{t+1}, a_{t+1}\i\}$. If $x_t$  and $\tau$ do not intersect, then
$\tau_1=\tau$, and we are done.
If $\mu(b_t, a_{t+1}\i) \ne 0$, then $x_t$  and $\tau$ intersect
only when $\mu(a_{t+1}\i,e) \ne 0$ for some $e \in 
E \backslash \{b_t, a_t\i, b_t\i,a_{t+1}, b_{t+1}, a_{t+1}\i \}$. In that case,
$\out(\tau_1)=\out(\tau) \cup \{(b_t\i, a_{t+1}) \}$, 
therefore $|\out(\tau_1)|\le 4\eg-2$.  So suppose $\mu(b_t, a_{t+1}\i) = 0$
as well.  Applying $h$ may create new outer branches only of one of the
following types:
$$(b_t, a_t\i), (b_t\i, a_{t+1}), (a_{t+1}\i, b_{t+1}\i).$$
Since $\mu(a_t,b_t)=0$, we have  $\mu_1(a_t,b_t)=0$.
If  any of $\mu(a_t\i, b_t\i),\mu(a_{t+1}, b_{t+1}),
\mu(b_{t+1},a_{t+1}\i)$ are 0, then they will be 0 with $\mu_1$ 
instead of $\mu$ and we are done. So let's assume they are all non-zero.
Let's look at the case $\eg \ge 4$,
 since the argument is easiest in this case. 
Because  $\tau$ is of
 reduced-length, 
one of the outer branches which does not intersect any of the simple closed
curves $x_t,a_t,b_t,a_{t+1},b_{t+1}$
must have zero measure, and this is going to stay zero in $\mu_1$.
This gives $n(\tau_1)<4\eg-1$. Now lets look at the case $\eg=3$, 
 and without loss of generality assume $t=1$. If $\mu(a_2\i,b_2\i)\ne 0$,
 then again one of the same type of outer branches must have 0 measure, 
and again we are done.
Therefore assume $\mu(a_2\i,b_2\i)= 0$. The assumptions force 
$$\out(\tau)=E(R) \backslash \{(a_1,b_1),(b_1\i,a_2),(a_2\i, b_2\i) \},$$
but this is not a reduced-length train-track. This takes care of the case 
$\eg=3$.  Now look at the case $\eg=2$.
Similar to the case of $\eg=3$, it follows that 
 $\mu(a_1,b_2\i)=\mu_1(a_1,b_2\i)=0$, and we are done.

(ii) There are bad curves but $\mu(a_t,b_t)=0$. The 
existence of bad curves means that $\mu(b_t,e)\ne 0$
 for some  $e \in E(R)\backslash  \{a_t\i,b_t\i, a_{t+1},
b_{t+1}, a_{t+1}\i\}$. In this case the train-track 
obtained by pushing the bad curves across $\partial R$
is collapsible to $\tau_1$. Notice that $(a_t,b_t) \notin \out(\tau_1)$
and  $\out(\tau_1) \backslash
\out(\tau)$ may only contain $(b_t\i,a_{t+1}), (a_{t+1},b_{t+1}),
(a_{t+1}\i, b_{t+1}\i)$.
 Therefore as in $(i)$, if $\eg \ge 4$ we are done. If $\eg=2$ or $3$ and
say $t=1$, then one of the branches in 
$$E(R) \backslash \{(a_1,b_1), (b_1\i, a_2),(a_2,b_2), (a_2\i, b_2\i) \}$$
must be missed by $\out(\tau)$ (since $\tau$ is of reduced-length)
and it will be missed by $\out(\tau_1)$ as well. 

(iii) There are bad curves and  $\mu(a_t,b_t)\ne 0$.
In this case after pushing the bad curves, we still have to push some
``bad pairs'' which come out near the edge $a_{t+1}\i$.
In this case
$$\out(\tau_1)\backslash \out(\tau) \subseteq \{ b_t\i, a_{t+1}\i)\}.$$
Since $|\out(\tau)| \le 4 \eg-3$, $|\out(\tau_1)| \le 4 \eg-2$
and we are done. $\spadesuit$

\proclaim {Theorem 5.4} The complexity of the word problem in $\eM^0_\eg$
is $O(|w|^2 \eg^2+ |w| \eg^2 \log \eg)$, where $|w|$ is the word length in the
set of generators (1.2).
In particular, for  $|w| \ge \log \eg$, the word problem has complexity
$O(|w|^2 \eg^2)$.
\endproclaim

{\demo{Proof}} Since the word problem in $\eM_2$ is quadratic in the word length, we need to prove the theorem for $\eg\ge 3$. (This is because there are
mapping classes in $\eM_2$ which fix all simple closed curves but are not 
the identity element, so our methods purely do not solve the word problem).
 Put the analog of 
each $\nu_i$, $i=1,\cdots,4$ given  in Theorem 2.4 in 
a reduced-length form. This takes $O(\eg)$ since $T(\nu_i)=O(\eg)$.
 Given the word $w=h_1 \cdots h_{|w|}$, apply each generator 
on the $\nu_i$, $i=1,\cdots,4$. After each application
put the resulting measured train-track in a reduced-length form.
This takes $O(\eg^2 \log ({\text{size}}))$. But the size grows by at most a factor of 3, therefore the total complexity is
$$O( \eg^2\log(\eg)+\eg^2\log(3\eg)+\cdots+\eg^2\log(3^{|w|-1}\eg))$$
which is $O(|w|^2 \eg^2+ |w| \eg^2 \log \eg)$. $\spadesuit$

\

\head Appendix: Turing Machine and Computational Complexity \endhead

A Turing Machine (see [Br] or [S], for example) is a hypothetical machine
consisting of an infinitely long
tape, a read/write head connected to a control mechanism. The tape is divided into infinitely many cells, each of which contains a symbol from a finite alphabet (the alphabet contains a special symbol for blank cell).
The cells are scanned one at a time using the read/write head, which can write a new symbol on the cell just read, move in either direction or not move at all. At any given time, the machine  is in one of the finitely many internal states. The behavior of

  the machine and a possible change of state depends on the current state, and the symbol read from the tape. 

Formally, let $X \subset Y$ be finite alphabets. A {\it{Turing Machine}}
 is a 
quadruple $(Q, \delta, q_0, q_F)$ where $Q$ is a finite set of states, 
$\delta$ is a function defined on a subset of $Q \times Y$ to 
$Q \times Y \times \{ L, R, 0 \}$ which is the {\it{state transition function, }}
$q_0 \in Q$ is the {\it{ start state,}} and $q_F \in Q$ is the {\it{halt state.}} The symbols $L,R,0$ should be interpreted as moving the head to the left, right, or no move at all, respectively. The set $X$ is the {\it{input alphabet}}.

Intuitively, any problem which is solvable by a finite instruction set
is solvable by a Turing Machine (see Church's Thesis say in [Br]). Therefore, we only describe a "program" for our solutions.

To define the complexity of an algorithm, there isn't a unique way. We 
have chosen the complexity to be the number of steps the Turing Machine 
takes to come up with the answer. 

To compute an upper bound for the complexity of a problem, we add up the number of steps needed for each sub-problem. They are all computed according to the following idea: 
To input a number of size $N$ into the machine takes $\log N$ steps.
The reason is one can write it in base 2,
with $O(\log_2 N)$ digits. Also, to add two numbers of size $\le N$ takes
$\log N$ steps as well.  Now one can devise a Turing machine to add the numbers in $O(\log N)$ steps which we leave as an exercise. 

From a theoretical point of view this definition (or any equivalent one
with respect to complexity) seems appropriate since a Turing Machine
is in a sense the most basic computer. 
In a Random Access Memory machine (say a typical PC), One assumes that it
takes a constant
time to add {\it{any}} two numbers. This assumption seems reasonable 
only when using machine-size numbers.

\

\head REFERENCES \endhead

\roster

\item"[B]" J. S. Birman, Braids, Links, and mapping class
groups,Princeton Univ. Press, 1974.

\item"[BKL]" J. S. Birman, K. H. Ko, and J. S. Lee, A new approach to the 
word problem and conjugacy problems in the braid groups,  To appear in Adv. Math. Available at xxx.lanl.gov as math.GT/9712211

\item"[BS]" J. S. Birman and C. Series, Algebraic linearity for an
automorphism of a surface group, J. Pure.  Appl. Alg. 52 (1988)
227-275.

\item"[Br]"  D. S. Bridges, Computability, a mathematical sketchbook,
GTM 146, Springer-Verlag, 1994.

\item"[E]" Epstein et al.,
 Word processing in groups,
 Jones and Bartlett Publishers,
 1992.

\item"[FLP]" A. Fathi, F. Laundenbach and V. Poenaru, Travaux de
Thurston sur les surfaces,
 Ast\'erisque 66-67,
 Soc. Math\'ematique de France,
 1979.

\item"[G]" F. A. Garside, The braid group and other groups, Oxford Quart. J. of Math., 20 (1969), 235-254.

\item"[HC]" H. Hamidi-Tehrani and Z.-H. Chen,
 Surface diffeomorphisms via train-tracks, Topology and its
Applications, 73 (1996), 141-167.

\item"[He]" G. Hemion,
 On classification of knots and 3-dimensional spaces, Oxford
University Press, New York, 1992.

\item"[Hu]" S. Humphries, Generators for the mapping class group, in
Lecture notes in Mathematics 722  (Springer, Berlin) 44-47.

\item"[J]" D. L. Johnson, Topics in the theory of group presentations,
London Math. Soc. Lecture note Ser. 42 (Cambridge Univ. Press, 1980).

\item"[KM]" R. Kleinberg and W. Menasco, Train tracks and zipping
sequences for pseudo-Anosov braids, preprint. Available at
http://math.Cornell.edu/$\sim$rdk.

\item"[MM1]" H. Masur and Y. Minsky, Geometry of the complex of curves
I: Hyperbolicity. Preprint. Available at
http://www.math.sunysb.edu/$\sim$yair.

\item"[MM2]"  H. Masur and Y. Minsky, Geometry of the complex of curves
II: hierarchal Structure. Preprint. Available at
http://www.math.sunysb.edu/$\sim$yair.

\item"[Mo1]" L.  Mosher, Classification of pseudo-Anosovs, in Low
dimensional topology and Kleinian groups, (1986) 13-75.

\item"[Mo2]" L. Mosher, Mapping class groups are automatic, Ann. of
Math. 142 (1995) No. 2, 303-384.

\item"[Mo3]" L. Mosher, A user's guide to the mapping class group:
Once punctured surfaces, In Geometric and computational perspectives
on infinite groups, DIMACS Ser. Discrete Math. Theoret. Comput. Sci.,
25 (1996) 101-174.

\item"[P]" R.C. Penner, The action of the mapping class group on curves
in surfaces, L'Enseigment Mathematique 30 (1984), 39-55.

\item "[PH]" R. C. Penner with J. L. Harer, Combinatorics of train
tracks, Ann. of Math. Stud., vol. 125, Princeton Univ. Press,
Princeton, NJ, 1995.

\item"[S]" A. Salomaa, Computation and automata, Cambridge University Press, 1985.

\endroster

\enddocument